\documentclass[11pt, reqno]{amsart} 
\usepackage[utf8]{inputenc}

\usepackage{amssymb, amsmath, amsthm}
\usepackage[bookmarks, bookmarksdepth=2, colorlinks=true, linkcolor=blue, citecolor=blue, urlcolor=blue]{hyperref}

\usepackage[alphabetic,lite]{amsrefs}
\usepackage{verbatim}
\usepackage{amscd}   
\usepackage[all]{xy} 
\usepackage{youngtab} 
\usepackage{young} 
\usepackage{ytableau}
\usepackage{tikz}
\usepackage{ mathrsfs }
\usepackage{cases}
\usepackage{array}
\usepackage{cellspace}
\usepackage{calligra,mathrsfs}
\usepackage{bm}
\usepackage{graphicx}
\usepackage{rank-2-roots}
\usepackage{float}
\usepackage{enumitem}

\newcommand{\defi}[1]{{\bf\upshape\sffamily #1}}
\DeclareMathOperator{\End}{End}

\DeclareMathOperator{\ShHom}{\mathscr{H}\text{\kern -3pt {\calligra\large om}}\,}

\renewcommand{\a}{\alpha}
\renewcommand{\b}{\beta}
\newcommand{\bw}{\bigwedge}
\renewcommand{\d}{\delta}
\newcommand{\D}{\mathcal{D}}

\DeclareMathOperator{\HH}{H}

\def\kk{{\mathbf k}}

\newcommand{\onto}{\twoheadrightarrow}
\newcommand{\oo}{\otimes}
\newcommand{\OO}{\mathcal{O}}
\newcommand{\pd}{\partial}

\newcommand{\s}{\sigma}

\renewcommand{\SS}{\mathcal{SS}}

\newcommand{\Ext}{\operatorname{Ext}}

\newcommand{\GL}{\operatorname{GL}}

\newcommand{\Spec}{\operatorname{Spec}}
\newcommand{\Sym}{\operatorname{Sym}}

\newcommand{\coker}{\operatorname{coker}}
\renewcommand{\det}{\operatorname{det}}

\renewcommand{\ker}{\operatorname{ker}}

\newcommand{\bb}[1]{\mathbb{#1}}

\newcommand{\mc}[1]{\mathcal{#1}}
\newcommand{\mf}[1]{\mathfrak{#1}}
\newcommand{\ol}[1]{\overline{#1}}
\newcommand{\op}[1]{\operatorname{#1}}

\def\PP{{\mathbf P}}
\def\lra{\longrightarrow}

\newtheorem{theorem}{Theorem}[section]
\newtheorem*{theorem*}{Theorem}
\newtheorem*{problem*}{Problem}
\newtheorem{lemma}[theorem]{Lemma}

\newtheorem{corollary}[theorem]{Corollary}
\newtheorem*{corollary*}{Corollary}

\newtheorem*{main-thm*}{Main Theorem}
\newtheorem*{linear-resolutions*}{Theorem on Linear Resolutions}
\newtheorem*{regularity-powers*}{Theorem on Regularity}
\newtheorem*{injectivity-Ext*}{Theorem on Injectivity of Maps of Ext Modules}
\newtheorem*{Kodaira*}{Kodaira Vanishing for Determinantal Thickenings}

\theoremstyle{definition}

\newtheorem*{definition*}{Definition}
\newtheorem{example}[theorem]{Example}

\theoremstyle{remark}
\newtheorem{remark}[theorem]{Remark}
\newtheorem*{remark*}{Remark}

\numberwithin{equation}{section}

\begin{document}

\title{The splitting type of jet bundles on the projective line}

\author{Annet Kyomuhangi}
\address{Department of Mathematics, Busitema University, P.O. Box 236, Tororo}
\email{annet.kyomuhangi@gmail.com}

\author{Emanuela Marangone}
\address{Dipartimento di Matematica, Università degli Studi di Genova, Via Dodecaneso, 35, 16146 Genova, Italy \newline
\indent CIMAT - Centro de Investigación en Matemáticas, Valenciana, 36023 Guanajuato, Mexico}
\email{emanuela.marangone@cimat.mx}

\author{Claudiu Raicu}
\address{Department of Mathematics, University of Notre Dame, 255 Hurley, Notre Dame, IN 46556\newline
\indent Institute of Mathematics ``Simion Stoilow'' of the Romanian Academy}
\email{craicu@nd.edu}

\author{Ethan Reed}
\address{Academy of Mathematics and Systems Science, Chinese Academy of Sciences, No. 55 Zhongguancun
East Road, Beijing, 100190, China.}
\email{ethan.reed@amss.ac.cn}


\date{\today}


\begin{abstract} 
Given a vector bundle $\mc{E}$ on a smooth variety $X$, one can associate to it its jet bundles, or bundles of principal parts. These bundles encode the Taylor expansions of local sections of $\mc{E}$ and arise naturally in the study of differential operators. When $X=\PP^1$ is the projective line, jet bundles decompose as direct sums of line bundles, making it natural to study their splitting types. The splitting type is well understood in characteristic zero, but exhibits subtler behavior in positive characteristic. Our goal is to give a complete answer in every positive characteristic, and in particular to provide a characterization of semistability for jet bundles on $\PP^1$. We also determine their extremal summands and give applications to differential operators and formulas for cohomology of line bundles on infinitesimal thickenings of the diagonal.
\end{abstract}

\maketitle

\section{Introduction}\label{sec:intro}

Given a vector bundle $\mc{E}$ on a smooth variety $X$, we consider the associated \defi{jet bundles} (or \defi{vector bundles of principal parts}) $\mc{P}^k(\mc{E})$, $k\geq 0$, that describe Taylor expansions up to order $k$ for germs of sections of $\mc{E}$. When $X=\PP^1$ is the projective line over some field $\kk$, jet bundles split as direct sums of line bundles, and it is natural to study their splitting type. In characteristic $p=0$, the situation is well-understood \cites{biswas,dr-som,per,maak08}: we have for instance that if $\mc{E}=\mc{O}_{\PP^1}(d)$ is a line bundle and if $d\geq k$ or $d\leq -1$, then $\mc{P}^k(\mc{E}) \simeq \mc{O}_{\PP^1}(d-k)^{\oplus(k+1)}$ is semistable. Understanding the splitting type in characteristic $p>0$ is more subtle, and a first step was taken in \cite{maak}. The goal of our paper is to provide a complete answer to this question -- since vector bundles on $\PP^1$ split as direct sums of line bundles, it suffices to consider the case $\mc{E}=\mc{O}_{\PP^1}(d)$.

\begin{theorem}\label{thm:sigma}
    Let $p=\op{char}(\kk)$ and consider $k\geq 0$ and $d\in\bb{Z}$. There exists an involution $\sigma=\sigma^{p,d}$ of the set $\{0,\cdots,k\}$ such that
    \[\mc{P}^k(\mc{O}_{\PP^1}(d)) \simeq \bigoplus_{i=0}^k \mc{O}_{\PP^1}(d-i-\sigma(i)).\]
\end{theorem}

Using the cycle decomposition of $\sigma$, we get that a transposition $(i,j)$ contributes two copies of $\mc{O}_{\PP^1}(d-i-j)$, while a fixed point $i$ contributes a single copy of $\mc{O}_{\PP^1}(d-2i)$. To make Theorem~\ref{thm:sigma} explicit, it remains to describe the involution $\sigma$. As we explain in Section~\ref{subsec:pparts}, one instance when $\sigma$ does not depend on the underlying characteristic $p$ occurs in the range $-1\leq d\leq k$. In this case we have
\begin{equation}\label{eq:d-1-to-k}
\sigma^{p,d}(i) = \begin{cases}
    d-i & \text{if }0\leq i\leq d, \\
    d+k+1-i & \text{if }d+1\leq i\leq k,
\end{cases}\quad\text{hence }\mc{P}^k(\OO_{\PP^1}(d)) \simeq \OO_{\PP^1}^{\oplus (d+1)} \oplus \OO_{\PP^1}(-k-1)^{\oplus (k-d)}.
\end{equation}
The remaining cases in characteristic $p=0$ correspond to semistable jet bundles, where we have 
\[\sigma^{0,d}(i) = k-i \quad\text{for }i=0,\cdots,k,\quad\text{when }d\leq -1\text{ or }d\geq k.\]
To describe the involution $\sigma$ in characteristic $p>0$ we define a family of integer-valued functions (depending on $p$) $\xi^s_r:\{0,\cdots,r-1\}\lra\bb{Z}$ for $s\in\bb{Z}$, $r\geq 1$, which are characterized by the following properties:
\begin{enumerate}
    \item\label{rule1} If $0\leq s\leq r$ then
    \begin{equation}\label{eq:basecase-xi}
    \xi^s_r(i) = \begin{cases}
        s-r & \text{if }0\leq i<s,\\
        s & \text{if }s\leq i<r.
    \end{cases}
    \end{equation}
    \item\label{rule2} If $r\leq q=p^e$ then $\xi^{s+q}_r=\xi^s_r$.
    \item\label{rule3} If $r<s\leq q=p^e$, we let $q'=p^{e-1}$ and write $r=r_0q'+r_1$, $i=i_0q'+i_1$, with $0\leq r_1,i_1<q'$. Then
    \[\xi^s_r(i) = \begin{cases}
        \xi^s_{r_1}(i_1) & \text{if }0\leq i_1<r_1,\\
        \xi^{s-r}_{q'-r_1}(i_1-r_1) & \text{if }r_1\leq i_1<q'.
    \end{cases}\]
\end{enumerate}

\begin{theorem}\label{thm:sigma-from-xi}
    If $p=\op{char}(\kk)>0$ then the involution $\sigma=\sigma^{p,d}$ in Theorem~\ref{thm:sigma} is given by
    \begin{equation}\label{eq:def-sigma-from-xi} 
    \sigma(i) = \xi^{d+1}_{k+1}(i)+(k-i)\quad\text{for }i=0,\cdots,k.
    \end{equation}
\end{theorem}

In the case $-1\leq d\leq k$, $\xi^{d+1}_{k+1}$ is determined by \eqref{eq:basecase-xi}, so the formula for $\sigma$ in Theorem~\ref{thm:sigma-from-xi} agrees with~\eqref{eq:d-1-to-k}. One immediate application of Theorem~\ref{thm:sigma-from-xi} is to characterize the semistability of $\mc{P}^k(\mc{O}_{\PP^1}(d))$, which is equivalent to the condition $\xi^{d+1}_{k+1}(i)=0$ for all $i$ (and is trivial when $k=0$, since $\mc{P}^0(\mc{O}_{\PP^1}(d))=\mc{O}_{\PP^1}(d)$).

\begin{theorem}\label{thm:semi-stability}
    In characteristic zero, $\mc{P}^k(\mc{O}_{\PP^1}(d))$ is semistable if and only if $d\leq -1$ or $d\geq k$.

    Suppose that $\op{char}(\kk)=p>0$ and $k\geq 1$. Let $q=p^e>k$ with $e$ minimal, and set $q'=p^{e-1}$. If
    \[k = k_0q'+k_1,\quad d\equiv d_0q'+d_1\ (\text{mod }q),\quad\text{ with }0\leq k_0,d_0<p,\quad 0\leq k_1,d_1<q',\]
    then $\mc{P}^k(\mc{O}_{\PP^1}(d))$ is semistable if and only if $d_0\geq k_0$ and $d_1\in\{k_1,q'-1\}$.
\end{theorem}

We illustrate the prior results with two examples.

\begin{example}[First principal parts]\label{ex:first-principal-parts}
Suppose that $p>0$ and $k=1$. The classical splitting of the first jet bundle (see for instance\cite{maak}*{Theorem~7.1}) is given by
\[
\mc{P}^1(\mc{O}_{\PP^1}(d))\simeq
\begin{cases}
\mc{O}_{\PP^1}(d)\oplus\mc{O}_{\PP^1}(d-2) & \text{if }p\mid d,\\
\mc{O}_{\PP^1}(d-1)^{\oplus2} & \text{if }p\nmid d.
\end{cases}
\]
In particular, $\mc{P}^1(\mc{O}_{\PP^1}(d))$ is semistable precisely when $p\nmid d$. To recover this using our results, note that rule \eqref{rule2} and \eqref{eq:def-sigma-from-xi} imply that $\sigma$ only depends on the residue class of $d\text{ mod }p$. If $d=0$, \eqref{eq:basecase-xi} implies $\xi^{1}_{2}(0)=-1$, $\xi^1_2(1)=1$, hence $\sigma(0)=0$, $\sigma(1)=1$ and we get corresponding summands $\mc{O}_{\PP^1}(d)$ and $\mc{O}_{\PP^1}(d-2)$ for $\mc{P}^1(\mc{O}_{\PP^1}(d))$. If $d=1$ then formula \eqref{eq:basecase-xi} yields $\xi^{2}_{2}(0)=\xi^2_2(1)=0$, hence $\sigma$ is the transposition $(0,1)$ and $\mc{P}^1(\mc{O}_{\PP^1}(d))$ is semistable. If $2\leq d\leq p-1$ then for $i=0,1$ we get
\[\xi^{d+1}_2(i) \overset{\eqref{rule3}}{=}\xi^{d-1}_1(0)\overset{\eqref{rule2}}{=}\xi^{1}_1(0) \overset{\eqref{rule1}}{=}0\]
hence $\mc{P}^1(\mc{O}_{\PP^1}(d))=\mc{O}_{\PP^1}(d-1)^{\oplus2}$ is again semistable.
\end{example}

\begin{example}\label{ex:higher-principal-parts}
Let $p=3$, $k=4$, and $d=6$. Setting $r=k+1=5$, $s=d+1=7$, we evaluate $\xi^7_5(i)$ using rule \eqref{rule3} with $q=9$, $q'=3$. We get $\xi^7_5(i)=\xi^7_2(i)=\xi^1_2(i)$ for $i=0,1$, $\xi^7_5(i)=\xi^7_2(i-3)=\xi^1_2(i-3)$ for $i=3,4$, and $\xi^7_5(2)=\xi^2_1(0)$. Combining this with calculations from Example~\ref{ex:first-principal-parts} we get
\[
\renewcommand{\arraystretch}{1.15}
\begin{array}{c|rrrrr}
i & 0 & 1 & 2 & 3 & 4\\ \hline
\xi_5^7(i) & -1 & 1 & 0 & -1 & 1\\
\sigma(i) & 3 & 4 & 2 & 0 & 1\\
\end{array}
\] 
Theorems~\ref{thm:sigma} and~\ref{thm:sigma-from-xi} then provide the splitting
\[
\mc{P}^4(\mc{O}_{\PP^1}(6))\simeq
\mc{O}_{\PP^1}(3)^{\oplus2}\oplus\mc{O}_{\PP^1}(2)\oplus\mc{O}_{\PP^1}(1)^{\oplus2}.
\]
Note that since $d\geq k$, the characteristic zero splitting is instead $\mc{P}^4(\mc{O}_{\PP^1}(6))\simeq\mc{O}_{\PP^1}(2)^{\oplus5}$. We will revisit this calculation in Example~\ref{ex:equivariant-principal-parts} where we compute the torus equivariant splitting.
\end{example}

It follows from Theorem~\ref{thm:sigma} that every summand $\mc{O}_{\PP^1}(m)$ of $\mc{P}^k(\mc{O}_{\PP^1}(d))$ satisfies $m\geq d-2k$. Moreover, $\mc{O}_{\PP^1}(d-2k)$ occurs as a summand if and only if $\sigma(k)=k$, and it does so with multiplicity one. To understand this condition better, we write $\omega_{\PP^1}\simeq\mc{O}_{\PP^1}(-2)$ for the canonical bundle on $\PP^1$. There is a canonical short exact sequence
\begin{equation}\label{eq:can-ses-Pk-Pk-1}
0\lra \omega_{\PP^1}^{\oo k}(d) \lra \mc{P}^k(\mc{O}_{\PP^1}(d)) \lra \mc{P}^{k-1}(\mc{O}_{\PP^1}(d))\lra 0,
\end{equation}
and since $\omega_{\PP^1}^{\oo k}(d)\simeq\mc{O}_{\PP^1}(d-2k)$, it follows that $\mc{O}_{\PP^1}(d-2k)$ is a summand in $\mc{P}^k(\mc{O}_{\PP^1}(d))$ if and only if the sequence \eqref{eq:can-ses-Pk-Pk-1} splits (in which case the splitting is unique).

\begin{theorem}\label{thm:split-ses-pparts}
    Suppose that $\op{char}(\kk)=p>0$, $k\geq 1$, and let $q=p^e>k$ with $e$ minimal. The sequence \eqref{eq:can-ses-Pk-Pk-1} is split exact if and only if $d\equiv k-1\ (\op{mod }q)$.
\end{theorem}

Note that if $k=1$ then $q=p$ and Theorem~\ref{thm:split-ses-pparts} agrees with Example~\ref{ex:first-principal-parts}. In characteristic zero \eqref{eq:can-ses-Pk-Pk-1} splits precisely when $d=k-1$, and the canonical splitting was constructed
in \cite{biswas}*{Corollary 2.3}. To go further, it is natural to study the extremal twists for the line bundle summands of each $\mc{P}^k(\mc{O}_{\PP^1}(d))$. We define 
\begin{equation}\label{eq:def-m-pl-min}
    m_-(k,d)=\min_i(a_i),\ m_+(k,d)=\max_i(a_i),\text{ where }\mc{P}^k(\mc{O}_{\PP^1}(d))\simeq\bigoplus_{i=0}^k\mc{O}_{\PP^1}(a_i).
\end{equation}
In the semistable case characterized by Theorem~\ref{thm:semi-stability} we get $m_-(k,d)=m_+(k,d)=d-k$, and as we explain next the semistable case controls the extremal twists in general. We define 
\begin{equation}\label{eq:def-SSk}
    \SS_k = \{d\in\bb{Z}:\mc{P}^k(\mc{O}_{\PP^1}(d))\text{ is semistable}\},
\end{equation}
\begin{equation}\label{eq:def-D+-}
    D_-(k,d)=\max\{i\in\SS_k,\ i\leq d\},\quad D_+(k,d)=\min\{i\in\SS_k,\ i\geq d\}.
\end{equation}

\begin{theorem}\label{thm:extr-summands-Pkd}
    For $k\geq 0$ and $d\in\bb{Z}$, the extremal twists of line bundle summands of $\mc{P}^k(\mc{O}_{\PP^1}(d))$ are
    \[m_-(k,d)=D_-(k,d)-k \text{ and }m_+(k,d)=D_+(k,d)-k.\]
\end{theorem}

\begin{example}\label{ex:extremal-twists}
We let $p=3$, $k=4$, and $d=6$, as in Example~\ref{ex:higher-principal-parts}. Theorem~\ref{thm:semi-stability}, with $q=9$, yields
\[
\SS_4=\{d\in\bb{Z}:d\equiv4,5,7,8\ (\op{mod }9)\},
\]
thus \eqref{eq:def-D+-} implies $D_-(4,6)=5$ and $D_+(4,6)=7$. Using Theorem~\ref{thm:extr-summands-Pkd} we obtain
\[
m_-(4,6)=1,\qquad m_+(4,6)=3,
\]
which agrees with the splitting in Example~\ref{ex:higher-principal-parts}. This is in contrast with the characteristic zero situation, when we have $6\in\SS_4$, so $D_-(4,6)=D_+(4,6)=6$ and $m_-(4,6)=m_+(4,6)=2$.
\end{example}

The interest in studying jet bundles comes from their relationship with differential operators. Following \cite{grothendieck}*{Section~IV.16}, \cite{stacksProject}*{\href{https://stacks.math.columbia.edu/tag/0G3P}{Section 0G3P}}, we consider the \defi{sheaf $\D^k_{d,e}$ of order $\leq k$ differential operators} from $\mc{O}_{\PP^1}(d)$ to $\mc{O}_{\PP^1}(e)$, and the corresponding space of \defi{global differential operators} $\op{Diff}^k(\mc{O}_{\PP^1}(d),\mc{O}_{\PP^1}(e))$. We have
\begin{equation}\label{eq:diff-Dkde}
\D^k_{d,e}=\ShHom\left(\mc{P}^k(\mc{O}_{\PP^1}(d)),\mc{O}_{\PP^1}(e)\right),\quad\op{Diff}^k(\mc{O}_{\PP^1}(d),\mc{O}_{\PP^1}(e))=H^0\left(\PP^1,\D^k_{d,e}\right).
\end{equation}
As a consequence of our earlier results, we obtain the following.

\begin{corollary}\label{cor:diffops}
    Let $m_-=m_-(k,d)$ and $m_+=m_+(k,d)$ as in Theorem~\ref{thm:extr-summands-Pkd}. We have
    \begin{enumerate}
        \item $\D^k_{d,e}$ is globally generated if and only if $e\geq m_+$.
        \item $\op{Diff}^k(\mc{O}_{\PP^1}(d),\mc{O}_{\PP^1}(e))\neq 0$ if and only if $e\geq m_-$.
        \item For every $d\in\bb{Z}$, $\D^k_{d,d}$ is globally generated and
        \[\dim_{\kk}\op{Diff}^k(\mc{O}_{\PP^1}(d),\mc{O}_{\PP^1}(d))=(k+1)^2.\]
    \end{enumerate}
\end{corollary}

When $q=p^e$, $\op{Diff}^{q-1}(\mc{O}_{\PP^1}(d),\mc{O}_{\PP^1}(d))$ is a finite dimensional algebra, whose structure in the case $d=0$ is analyzed in \cite{smith}*{Section~3}. The techniques in \cite{smith} apply for a general $d$, and imply an identification
\[\op{Diff}^{q-1}(\mc{O}_{\PP^1}(d),\mc{O}_{\PP^1}(d)) = \End\left(F^q_*\mc{O}_{\PP^1}(d)\right),\]
where $F^q$ denotes the $q$-th power Frobenius on $\PP^1$. Using the well-known splitting of $F^q_*\mc{O}_{\PP^1}(d)$, one can then show that $\op{Diff}^{q-1}(\mc{O}_{\PP^1}(d),\mc{O}_{\PP^1}(d))$ is semisimple (isomorphic to the algebra of $q\times q$ matrices) precisely when $d\equiv -1\ (\op{mod }q)$, or equivalently, when $\mc{P}^{q-1}(\mc{O}_{\PP^1}(d))$ is semistable.

Our original motivation for studying jet bundles on $\PP^1$ came from some important connections to cohomology of line bundles on the incidence correspondence, the Green--Han--Monsky ring, and Lefschetz properties (see \cites{KMRR0,KMRR3,KMRR1,KMRR2}). Given a projective space $\PP$ with coordinates $[x_1:\cdots:x_n]$, and its dual $\PP^{\vee}$ with coordinates $[y_1:\cdots:y_n]$, the \defi{incidence correspondence}
\[ X = \{(p,H) : p\in H\} \subset \PP \times \PP^{\vee}\]
is the hypersurface cut out by the bilinear form
\begin{equation}\label{eq:def-omega-nvars}
\omega = x_1y_1+\cdots+x_ny_n.
\end{equation}
We can consider more generally the infinitesimal thickenings
$X^{(r)}\subset \PP\times\PP^{\vee}$, defined by $\omega^r=0$, and corresponding line bundles on $X^{(r)}$ obtained by restriction from $\PP \times \PP^{\vee}$. We write
\[ \mc{O}_{X^{(r)}}(a,b) = \mc{O}_{\PP \times \PP^{\vee}}(a,b)_{|_{X^{(r)}}}\quad\text{ for }(a,b)\in\bb{Z}^2.\]
The main results of \cite{KMRR1} emerge from understanding the cohomology of $\mc{O}_{X^{(r)}}(a,b)$ with $r\leq p$. When $n=2$, $\PP=\PP^1$ and there is a natural identification $\PP^{\vee}=\PP^1$ which makes $X$ the diagonal in $\PP^1\times\PP^1$. One gets an immediate relationship between the line bundles $\mc{O}_{X^{(r)}}(a,b)$ and the jet bundles on $\PP^1$, and our computation of the splitting type provides explicit formulas for the cohomology of $\mc{O}_{X^{(r)}}(a,b)$. We discuss this in Section~\ref{sec:coh-thickenings}, and record here the following consequence.

\begin{theorem}\label{thm:coho-consequences}
    If $n=2$ and $r\geq 1$, then $\mc{O}_{X^{(r)}}(a,b)$ has vanishing cohomology if and only if 
    \[a+b=r-2\quad\text{ and }\quad\mc{P}^{r-1}(\mc{O}_{\PP^1}(b))\text{ is semistable}.\]
    In characteristic $p>0$, if we write $h^i$ for the dimension of the $i$-th cohomology group, then
    \[h^0\left(X^{(r)},\mc{O}_{X^{(r)}}(a,b)\right)=h^1\left(X^{(r)},\mc{O}_{X^{(r)}}(a,b)\right)=\frac{1}{2}\sum_{i=0}^{r-1}|\xi^{b+1}_r(i)|\quad\text{ whenever }a+b=r-2.\]
\end{theorem}

\begin{example}\label{ex:cohomology-thickening}
Let $p=3$ and $r=5$, and consider $(a,b)=(-3,6)$, for which $a+b=r-2$. As discussed in Example~\ref{ex:higher-principal-parts}, $\mc{P}^4(\mc{O}_{\PP^1}(6))$ is not semistable. The computation of $\xi^7_5(i)$ and Theorem~\ref{thm:coho-consequences} then imply
\[h^0\left(X^{(5)},\mc{O}_{X^{(5)}}(-3,6)\right)=h^1\left(X^{(5)},\mc{O}_{X^{(5)}}(-3,6)\right)=2.\]
By contrast, in characteristic zero $\mc{P}^4(\mc{O}_{\PP^1}(6))$ is semistable and $\mc{O}_{X^{(5)}}(-3,6)$ has vanishing cohomology.
\end{example}

The main technical aspect of our approach to understanding the jet bundles $\mc{P}^k(\mc{O}_{\PP^1}(d))$ comes from an analysis of related vector bundles $\mc{F}^d_r$ introduced in \eqref{eq:def-seq-Fdr}. When $d\geq r-1$, we have that $\mc{F}^d_r\simeq\D^{r-1}_{d,0}$. In Theorem~\ref{thm:recursive-Fdr} we give a recursive description of the splitting type of $\mc{F}^d_r$, keeping track of the equivariant structure with respect to the natural action of the $2$-dimensional torus. Without much additional effort, the remaining statements follow as applications of this result.

\medskip

\noindent{\bf Organization.} In Section~\ref{sec:prelim} we collect a number of elementary observations. In Section~\ref{sec:pparts-P1} we provide the recursive description for the splitting type of the bundles $\mc{F}^d_r$. In Section~\ref{sec:split-jet} we explain how this leads to the main results about jet bundles on $\PP^1$. Section~\ref{sec:extremal-summands} is concerned with extremal twists for line bundle summands of $\mc{P}^k(\mc{O}_{\PP^1}(d))$, and consequences for differential operators. We end in Section~\ref{sec:coh-thickenings} with a discussion of applications to cohomology on infinitesimal thickenings.

\section{Preliminaries}\label{sec:prelim}

Throughout this paper $\kk$ is an algebraically closed field of characteristic $p$, $U=\kk^2$ is a $2$-dimensional vector space, and $\PP^1=\bb{P}U=\op{Proj}(\Sym U)$ is the projective line. It comes with a tautological short exact sequence
\begin{equation}\label{eq:ses-on-P1}
0 \lra \mc{R} \lra U \oo \mc{O}_{\PP^1} \lra \mc{O}_{\PP^1}(1) \lra 0,
\end{equation}
where $\mc{R}=\Omega^1_{\PP^1}(1) = \bw^2 U\oo \mc{O}_{\PP^1}(-1)$ is the universal subsheaf and $\mc{O}_{\PP^1}(1)$ is the universal quotient sheaf. The natural action of $\GL(U)$ on $U$ induces an action on $\PP^1$ and a canonical linearization of $\mc{O}_{\PP^1}(1)$ such that \eqref{eq:ses-on-P1} is an exact sequence of $\GL(U)$-equivariant sheaves.

\subsection{Some $\GL$-equivariant bundles on $\PP^1$}
\label{subsec:def-Fdr}

For $d\geq 0$, the sequence \eqref{eq:ses-on-P1} induces a short exact sequence
\[ 0 \lra \mc{R}^{\oo d} \lra D^dU \oo \mc{O}_{\PP^1} \overset{\Delta_d}{\lra} D^{d-1}U\oo\mc{O}_{\PP^1}(1) \lra 0,
\]
where $D^d$ is the \defi{divided power} functor, and $\Delta_d$ is induced by the comultiplication map $D^dU\lra D^{d-1}U\oo U$. We define a family of $\GL(U)$-equivariant locally free sheaves $\mc{F}^d_r$ for $r\geq 0$ via
\begin{equation}\label{eq:def-seq-Fdr} 0 \lra \mc{F}^d_r \lra D^dU \oo \mc{O}_{\PP^1} \overset{\Delta}{\lra} D^{d-r}U \oo \mc{O}_{\PP^1}(r) \lra 0,
\end{equation}
where $\Delta = (\Delta_{d-r+1}(r-1)) \circ \cdots \circ (\Delta_{d-1}(1)) \circ \Delta_d$. We have that $\Delta=0$ if $d<r$, hence
\begin{equation}\label{eq:Fdr-for-d<r} 
 \mc{F}^d_r = D^dU \oo \mc{O}_{\PP^1}\quad\text{when }d<r,
\end{equation}
which is why we will usually assume that $d\geq r$. In this case we have by construction that $\mc{F}^d_r$ has a filtration with composition factors
\[ \mc{R}^{\oo d},\ \mc{R}^{\oo(d-1)}(1),\cdots,\ \mc{R}^{\oo(d-r+1)}(r-1)\]
and therefore
\begin{equation}\label{eq:rank-det-Fdr}
\op{rank}(\mc{F}^d_r) = r \quad\text{ and }\quad \det(\mc{F}^d_r) = \bw^r\mc{F}^d_r = \left(\bw^2 U\right)^{\oo\left(dr-{r\choose 2}\right)} \oo \mc{O}_{\PP^1}(-r(d-r+1)).
\end{equation}
In particular, $\mc{F}^d_0=0$. Moreover, by the snake lemma we have natural short exact sequences
\begin{equation}\label{eq:ses-Fdrs}
    0 \lra \mc{F}^d_r \lra \mc{F}^d_s \lra \mc{F}^{d-r}_{s-r}(r) \lra 0\qquad\text{ for all }s\geq r.
\end{equation}
Note that when $s>d$, \eqref{eq:ses-Fdrs} reduces to \eqref{eq:def-seq-Fdr}.

\begin{lemma}\label{lem:H01-Fdr}
    For $d_1,d_2\geq r-1$ and $i=0,1$, there exist isomorphisms
\[H^i\left(\PP^1,\mc{F}^{d_1}_r(d_2-r)\right) \simeq H^{1-i}\left(\PP^1,\mc{F}^{d_2}_r(d_1-r)\right)^{\vee}.\]    
\end{lemma}

\begin{proof}
     Using \eqref{eq:def-seq-Fdr} and the long exact sequence in cohomology, we get an exact sequence
 \begin{equation}\label{eq:les-coh-Fdr}
0\to H^0\left(\PP^1,\mc{F}^{d_1}_r(d_2-r)\right) \to D^{d_1} U \oo \Sym^{d_2-r}U \overset{\Delta}{\to} D^{d_1-r}U \oo \Sym^{d_2} U\to H^1\left(\PP^1,\mc{F}^{d_1}_r(d_2-r)\right)\to 0.
 \end{equation}
Using $(\Sym^i U)^{\vee} = D^i(U^{\vee})$, $(D^i U)^{\vee} = \Sym^i(U^{\vee})$, and using the isomorphism $U^{\vee} \simeq \kk^2 \simeq U$, we obtain by dualizing \eqref{eq:les-coh-Fdr} the exact sequence
\[0\to H^1\left(\PP^1,\mc{F}^{d_1}_r(d_2-r)\right)^{\vee} \to D^{d_2} U \oo \Sym^{d_1-r}U \overset{\Delta}{\to} D^{d_2-r}U \oo \Sym^{d_1} U\to H^0\left(\PP^1,\mc{F}^{d_1}_r(d_2-r)\right)^{\vee}\to 0.
\]
Comparing this with the exact sequence obtained from \eqref{eq:les-coh-Fdr} by exchanging $d_1$ with $d_2$, we obtain the desired isomorphisms.
\end{proof}

In the proof of Theorem~\ref{thm:recursive-Fdr} we will also need the \defi{Frobenius power functor} $F^q$ and the \defi{$q$-truncated divided power functors} $T_qD^d$, where $p > 0$ and $q = p^e$ for $e\geq 0$. Here $F^qV\subseteq\Sym^q V$ denotes the linear span of the powers $v^q$ with $v\in V$, and
\[ T_qD^d V = \ker\left(D^d V \lra D^{d-q}V \oo F^q V\right).\]
For a line bundle $\mc{L}$, we have $T_qD^d\mc{L}=\mc{L}^d$ when $0\leq d<q$ and $T_qD^d\mc{L}=0$ otherwise.

\subsection{Torus equivariant splitting}
\label{subsec:torus-splitting}

By Grothendieck's theorem, each vector bundle $\mc{E}$ on $\PP^1$ splits as a direct sum of line bundles. When present, we would like to take advantage of the $\GL(U)$-equivariant structure of $\mc{E}$, but the splitting cannot usually be taken to be $\GL(U)$-equivariant. We therefore fix a basis of $U$, make the identification $U=\kk^2$, and consider the $2$-dimensional torus $T=\bb{G}_m \times \bb{G}_m \subset\GL(U)$ acting coordinate-wise on~$U$. By \cite{kumar}, there exists a decomposition
\[\mc{E} = \mc{L}_1 \oplus \cdots \oplus \mc{L}_r,\]
where each $\mc{L}_i$ is a $T$-equivariant line bundle (note that although the statement in \cite{kumar} is phrased over $\bb{C}$, the proof only uses the fact that irreducible $T$-representations are $1$-dimensional, and taking $T$-invariants is an exact functor). Each $T$-equivariant line bundle on $\PP^1$ is of the form
\[L_{u,v}\oo\mc{O}_{\PP^1}(i),\text{ where }u,v,i\in\bb{Z}
\]
and $L_{u,v}=\kk\cdot\chi^{u,v}$ is the $1$-dimensional $T$-representation where the action is given by
\begin{equation}\label{eq:act-on-chiuv}
 (t_1,t_2) * \chi^{u,v} = t_1^{u}t_2^{v}\chi^{u,v}\text{ for }t_1,t_2\in\bb{G}_m.
\end{equation}
For the standard representation $U$, we have the $T$-equivariant decomposition $U=L_{1,0}\oplus L_{0,1}$. The determinant of $U$ is then $\det U=\bigwedge^2 U \simeq L_{1,1}$. More generally, for $a\geq b$ we consider the $T$-representation
\begin{equation}\label{eq:def-Hab} 
H_{a,b} = \bigoplus_{i=b}^a L_{a+b-i,i},
\end{equation}
and observe that $H_{a,b}$ is isomorphic to $(\det U)^{\oo b}\oo\Sym^{a-b}U$ viewed as a $T$-representation. In particular $U=H_{1,0}$, $\Sym^d U = H_{d,0}$, and $\det(U)=H_{1,1}$. We make the convention $H_{a,b}=0$ if $a<b$.

For $p=\op{char}(\kk)>0$ and $q=p^e$, we consider the Frobenius action on $T$-representations, and more generally on $T$-equivariant sheaves on $\PP^1$, which is given on line bundles by
\[ F^q\left(L_{u,v}\oo\mc{O}_{\PP^1}(i)\right) = L_{qu,qv}\oo\mc{O}_{\PP^1}(qi).\]
For a $T$-representation $W$ and $i\in\bb{Z}$, we will often write $W(i)$ instead of $W\oo\mc{O}_{\PP^1}(i)$. 

\subsection{Principal parts on $\PP^1$}
\label{subsec:pparts}

Let $k\geq 0$, and let $\mc{I}$ be the ideal sheaf of the diagonal in
$\PP^1\times\PP^1$. There is a natural $\GL(U)$-equivariant identification
$\mc{I}\simeq\det U\oo\OO_{\PP^1\times\PP^1}(-1,-1)$.
Write $\PP^{(k+1)}$ for the closed subscheme defined by $\mc{I}^{k+1}$, so that
\begin{equation}\label{eq:ses-pparts}
0\lra (\det U)^{k+1}\oo\OO_{\PP^1\times\PP^1}(-k-1,-k-1)
\lra\OO_{\PP^1\times\PP^1}\lra\OO_{\PP^{(k+1)}}\lra 0.
\end{equation}
Let $\pi_1,\pi_2$ denote the projections $\PP^1\times\PP^1\lra\PP^1$
and their restrictions to $\PP^{(k+1)}$. For $d\in\bb{Z}$, the
\defi{bundle of principal parts of order $k$} of $\OO_{\PP^1}(d)$ is
\[
\mc{P}^k(\OO_{\PP^1}(d))
=\pi_{1*}\pi_2^*\OO_{\PP^1}(d)
=\pi_{1*}\OO_{\PP^{(k+1)}}(0,d).
\]
Twisting \eqref{eq:ses-pparts} by $\OO_{\PP^1\times\PP^1}(0,d)$
and pushing forward along $\pi_1$ gives the following three cases.

\smallskip

\noindent{\bf Case 1.} If $d\geq k$ then we get a short exact sequence
\begin{equation}\label{eq:ses-Pkm-case1}
0 \lra (\det U)^{k+1} \oo \Sym^{d-k-1} U \oo \OO_{\PP^1}(-k-1) \lra \Sym^d U \oo \OO_{\PP^1} \lra \mc{P}^k(\OO_{\PP^1}(d)) \lra 0 
\end{equation}
Dualizing, tensoring by $(\det U)^d$, and using the identification $(\Sym^d U)^{\vee} \oo (\det U)^d = D^d U$, we get an exact sequence
\[ 0 \lra \left(\mc{P}^k(\OO_{\PP^1}(d))\right)^{\vee} \oo (\det U)^d \lra D^d U \oo \OO_{\PP^1} \lra D^{d-k-1}U \oo \OO_{\PP^1}(k+1) \lra 0.\]
Taking $r=k+1$ and comparing with \eqref{eq:def-seq-Fdr} we conclude that
\begin{equation}\label{eq:Pkm-vs-Fdr-case1} \mc{P}^k(\OO_{\PP^1}(d)) = \left(\mc{F}^d_{k+1}\right)^{\vee} \oo (\det U)^d = \left(\mc{F}^d_{k+1}\right)^{\vee} \oo L_{d,d} 
\end{equation}
so the ($T$-equivariant) splitting of $\mc{P}^k(\OO_{\PP^1}(d))$ is determined by the calculations in Section~\ref{sec:pparts-P1}.

\smallskip

\noindent{\bf Case 2.} If $-1\leq d\leq k$ then we get a short exact sequence
\[ 0 \lra \Sym^d U \oo \OO_{\PP^1} \lra \mc{P}^k(\OO_{\PP^1}(d)) \lra (\det U)^k \oo (\Sym^{k-1-d}U)^{\vee} \oo \OO_{\PP^1}(-k-1) \lra 0\]
Since $\Ext^1(\OO_{\PP^1}(-k-1),\OO_{\PP^1})=\HH^1({\PP^1},\OO_{\PP^1}(k+1))=0$, the above sequence splits and we get
\[ 
\begin{aligned}
    \mc{P}^k(\OO_{\PP^1}(d)) &\simeq \left(\bigoplus_{i=0}^d L_{i,d-i}(0) \right)  \oplus \left(\bigoplus_{i=0}^{k-1-d} L_{d+1+i,k-i}(-k-1) \right) \\
    & \simeq \OO_{\PP^1}^{\oplus (d+1)} \oplus \OO_{\PP^1}(-k-1)^{\oplus (k-d)}
\end{aligned}
\]
Here the second isomorphism is nonequivariant.

\smallskip

\noindent{\bf Case 3.} If $d\leq -2$ then we get a short exact sequence
\[
\begin{aligned}
0 \lra \mc{P}^k(\OO_{\PP^1}(d)) &\lra (\det U)^k \oo (\Sym^{k-1-d}U)^{\vee} \oo \OO_{\PP^1}(-k-1) \lra \\
&\lra (\det U \oo \Sym^{-d-2}U)^{\vee} \oo \OO_{\PP^1} \lra 0.
\end{aligned}
\]
Tensoring by $(\det U)^{-1-d} \oo \OO_{\PP^1}(k+1)$ and passing to divided powers as in {\bf Case 1} we get a short exact sequence
\[0 \lra \mc{P}^k(\OO_{\PP^1}(d)) \oo (\det U)^{-1-d} \oo \OO_{\PP^1}(k+1) \lra D^{k-1-d}U \oo \OO_{\PP^1} \lra D^{-2-d}U \oo \OO_{\PP^1}(k+1) \lra 0 \]
Comparing with \eqref{eq:def-seq-Fdr} we conclude that
\begin{equation}\label{eq:PkOd-neg-d} \mc{P}^k(\OO_{\PP^1}(d)) = \mc{F}^{k-1-d}_{k+1} \oo (\det U)^{1+d} \oo \OO_{\PP^1}(-k-1)\simeq \mc{F}^{k-1-d}_{k+1} \oo L_{1+d,1+d}(-k-1)
\end{equation}
and the ($T$-equivariant) splitting of $\mc{P}^k(\OO_{\PP^1}(d))$ is again determined by that of the bundles~$\mc{F}^d_r$.

\begin{remark}\label{rem:principal-parts}
    When $d\geq k$ and $\op{char}(\kk)=0$, it follows from \cite{dr-som}*{Proposition~2.2}, \cite{maak}*{Proposition~6.3}, or \cite{maak08}*{Theorem~1.1} that 
    \[\mc{P}^{k}(\mc{O}_{\PP^1}(d))\simeq\mc{O}_{\PP^1}(d-k)^{\oplus (k+1)}.\]
    Keeping track of the equivariant structure in \eqref{eq:ses-Pkm-case1}, we in fact get
    \begin{equation}\label{eq:PkOd-char0-semistable}
    \mc{P}^{k}(\mc{O}_{\PP^1}(d))=\Sym^k U \oo \mc{O}_{\PP^1}(d-k)=H_{k,0}(d-k),
    \end{equation}
    where the first identification is $\GL(U)$-equivariant, and the second one is $T$-equivariant. Setting $r=k+1$, it follows from \eqref{eq:Pkm-vs-Fdr-case1} that
\begin{equation}\label{eq:char-0-GLU-Fdr}
\mc{F}^d_r = (\det U)^d\oo(\Sym^{r-1}U)^{\vee}\oo\mc{O}_{\PP^1}(-(d-r+1))=H_{d,d-r+1}(r-1-d)\quad\text{ for }d\geq r-1.
\end{equation}
Combining \eqref{eq:PkOd-neg-d} and \eqref{eq:char-0-GLU-Fdr}, it follows that \eqref{eq:PkOd-char0-semistable} also holds when $d\leq -2$.
\end{remark}

\subsection{A local cohomology model}
\label{subsec:addl}

Let $x_1,x_2$ be the basis of $U$ fixed in Section~\ref{subsec:torus-splitting},
and let $y_1,y_2$ be the dual basis of $U^\vee$. Set
\[
S=\kk[x_1,x_2]=\Sym U,\qquad R=S[y_1,y_2],\qquad
M=\bw^2U^\vee\oo H^2_{(y_1,y_2)}(R).
\]
We give $M$ the bigrading
\begin{equation}\label{eq:explicit-basis-Mde}
M_{d,e}=D^dU\oo\Sym^eU
\simeq\bw^2U^\vee\oo
\bigoplus_{\substack{d_1,d_2,e_1,e_2\geq0\\d_1+d_2=d,\ e_1+e_2=e}}
\kk\cdot\frac{x_1^{e_1}x_2^{e_2}}{y_1^{1+d_1}y_2^{1+d_2}}
\qquad(d,e\geq0).
\end{equation}
A Laurent monomial represents zero whenever its exponent of $y_1$ or $y_2$
is nonnegative. From now on we suppress the factor $\bw^2U^\vee$ in the
notation, but retain its contribution to the $T$-action. In particular,
\begin{equation}\label{eq:forget-bw2Udual}
\kk\cdot\frac{1}{y_1^{1+d_1}y_2^{1+d_2}}\simeq L_{d_1,d_2}.
\end{equation}
Multiplication by
\begin{equation}\label{eq:def-omega}
\omega=x_1y_1+x_2y_2
\end{equation}
realizes the map induced by comultiplication on divided powers and
multiplication on symmetric powers:
\begin{equation}\label{eq:Delta-de-is-mult-omega}
M_{d,e}\overset{\cdot\omega}{\lra}M_{d-1,e+1}.
\end{equation}

\section{Splitting the bundles $\mc{F}^d_r$}\label{sec:pparts-P1}

The goal of this section is to describe recursively the $T$-equivariant splitting of $\mc{F}^d_r$ in characteristic $p>0$. 

\begin{theorem}\label{thm:recursive-Fdr}
    Suppose that $\op{char}(\kk)=p>0$, $d\geq0$, and $0\leq r\leq q=p^e$, where $e\geq0$. We have the following $T$-equivariant isomorphisms:
    \begin{enumerate}
        \item If $d<r$ then 
        \[ \mc{F}^d_r = H_{d,0}\oo\mc{O}_{\PP^1}.\]
        \item If $d\geq q+r-1$ then
        \[ \mc{F}^d_r = L_{q,q} \oo \mc{F}^{d-q}_r(-q).\]
        \item If $q-1\leq d \leq q+r-1$ then
        \[ \mc{F}^d_r = H_{d,q}(-q) \bigoplus H_{q-1,d-r+1}(r-q).\]
        \item If $e\geq1$, set $q'=p^{e-1}$. If $q'\leq r\leq d\leq q-1$ and $1\leq a<p$ are such that $aq'\leq r<(a+1)q'$, then
        \[\mc{F}^d_r = F^{q'}(H_{a,0}) \oo \mc{F}^{d-aq'}_{r-aq'} \bigoplus F^{q'}(H_{a-1,0}) \oo 
        \mc{F}^{d-r+q'}_{(a+1)q'-r}(r-aq').\]
    \end{enumerate}
\end{theorem}

The proof of Theorem~\ref{thm:recursive-Fdr} will occupy the rest of the section. Cases (1)--(3) are somewhat straightforward, but case (4) seems to require a more delicate analysis, and the crucial ingredient is to understand the action of differential operators on the sheaves $\mc{F}^d_r$, and their corresponding graded modules of sections $F^d_r$. Before going into the details of the proof, we point out the following reformulation of part (4).

\begin{corollary}\label{cor:rec-Fdr-sometimes}
 Suppose that $q=p^e$, $e\geq1$, and $q'=p^{e-1}\leq r\leq q$.
 \begin{enumerate}
  \item[(i)] If $1\leq a<p$ and $aq'\leq r\leq d < (a+1)q'$ then
        \[\mc{F}^d_r = F^{q'}(H_{a,0}) \oo \mc{F}^{d-aq'}_{r-aq'} \bigoplus F^{q'}(H_{a-1,0}) \oo \left[H_{d-r+q',q'}(r-(a+1)q') \bigoplus H_{q'-1,d-aq'+1}\oo\mc{O}_{\PP^1}\right] \] 
  \item[(ii)] In particular, if $p=2$ and $q'\leq r\leq d<q$ then 
          \[\mc{F}^d_r = F^{q'}U \oo \mc{F}^{d-q'}_{r-q'} \bigoplus H_{d-r+q',q'}(r-q) \bigoplus H_{q'-1,d-q'+1}\oo\mc{O}_{\PP^1} \] 
 \end{enumerate}
\end{corollary}

\begin{proof}
 Part (ii) is the special case of (i) when $p=2$ and $a=1$, so it suffices to consider case (i). Using Theorem~\ref{thm:recursive-Fdr}(4), the conclusion reduces to showing that
 \[ \mc{F}^{d-r+q'}_{(a+1)q'-r} = H_{d-r+q',q'}(-q') \bigoplus H_{q'-1,d-aq'+1}(aq'-r).\]
 If we set $d'=d-r+q'$, $r'=(a+1)q'-r$, then the hypothesis $aq'\leq r\leq d\leq (a+1)q'-1$ implies $q'-1\leq d'\leq q'+r'-1$, so we can apply Theorem~\ref{thm:recursive-Fdr}(3) to get the desired description of $\mc{F}^{d-r+q'}_{(a+1)q'-r} =\mc{F}^{d'}_{r'}$.
\end{proof}

We next introduce some notation to be used in the rest of the section. We let $S=\kk[x_1,x_2] = \Sym_{\kk}U$, write $S(r)$ for the rank-one graded module generated in degree $-r$, and define the graded $S$-module $F^d_r$ via
\[ F^d_r = \ker\left(D^d U \oo S \lra D^{d-r}U\oo S(r)\right),\]
induced by the comultiplication $D^d U \lra D^{d-r}U \oo D^r U$, followed by the natural map $D^r U \lra \Sym^r U$ and the multiplication $\Sym^r U \oo S \lra S(r)$. Notice that $\mc{F}^d_r$ is the sheaf associated to the graded module $F^d_r$. Using the notation in Section~\ref{subsec:addl}, we get that
\[ D^d U \oo S = M_{d,\bullet} = \bigoplus_{e\geq 0}M_{d,e}.\]
Moreover, for each $e$ we get an exact sequence
\[ 0 \lra H^0\left(\PP^1,\mc{F}^d_r(e)\right) \lra M_{d,e} \overset{\cdot\omega^r}{\lra} M_{d-r,e+r} \]
Using the explicit description \eqref{eq:explicit-basis-Mde}, we get that multiplication by $y_i$ induces an $S$-linear map $M_{d,\bullet}\lra M_{d-1,\bullet}$ which commutes with multiplication by $\omega$, hence we get induced maps
\begin{equation}\label{eq:mult-yi-Fdr}
 y_i : F^d_r \lra F^{d-1}_r\quad\text{ and }\quad y_i : \mc{F}^d_r \lra \mc{F}^{d-1}_r.    
\end{equation}

\subsection{Local charts for the sheaves $\mc{F}^d_r$.} It will be useful to have a local description of the sheaves $\mc{F}^d_r$, and because of symmetry we will only focus on the basic open affine $D_+(x_1)$, leaving the analysis in the chart $D_+(x_2)$ to the interested reader (see also \cite{maak}). We write $t=x_2/x_1$, $A=\kk[t]$ so that $D_+(x_1)=\Spec(A)$, and consider the identifications
\begin{equation}\label{eq:A-basis-Ddu}D^dU\oo\mc{O}_{\PP^1}(D_+(x_1)) = \op{Span}_A\left\{\frac{1}{y_1^{1+d_1}y_2^{1+d_2}}:d_1+d_2=d\right\},
\end{equation}
\[
D^{d-r}U\oo\mc{O}_{\PP^1}(r)(D_+(x_1)) = \op{Span}_A\left\{\frac{x_1^r}{y_1^{1+d_1}y_2^{1+d_2}}:d_1+d_2=d-r\right\}.\]

\begin{lemma}\label{lem:local-basis-Fdr}
    If $d\geq r$ then the following sections give an $A$-basis of $\mc{F}^d_r(D_+(x_1))$:
    \[ f_i = \sum_{j=i}^{d} {-r \choose j-i}\cdot t^{j-i} \cdot\frac{1}{y_1^{j+1}y_2^{d+1-j}},\quad\text{ for }i=0,\cdots,r-1.\]
\end{lemma}

\begin{proof} Relative to the basis of $D^dU\oo A$ from \eqref{eq:A-basis-Ddu}, the $f_i$'s are described by the columns of the matrix
\[
\begin{bmatrix}
   1 & 0 & 0 & \cdots & 0 \\
   {-r\choose 1}\cdot t & 1 & 0 & \cdots & 0 \\
   {-r\choose 2}\cdot t^2 & {-r\choose 1}\cdot t & 1 & \cdots & 0 \\
   \vdots & \vdots & \vdots & \ddots & \vdots \\
   \vdots & \vdots & \vdots & \ddots & 1 \\
   \vdots & \vdots & \vdots & \ddots & \vdots \\
  {-r\choose d}\cdot t^d & {-r\choose d-1}\cdot t^{d-1} & {-r\choose d-2}\cdot t^{d-2} & \cdots & {-r\choose d-r+1}\cdot t^{d-r+1} \\
\end{bmatrix}
\]
Since the submatrix formed by the first $r$ rows has determinant $1$, the columns span a rank $r$ direct summand of $D^dU\oo A$. To show that this agrees with $\mc{F}^d_r(D_+(x_1))$, it then suffices to check that each $f_i$ is annihilated by $\omega^r = (y_1+ty_2)^r\cdot x_1^r$. To that end, we consider the identity
\[ 1 = (1+z)^{-r}\cdot (1+z)^r = \left(\sum_{u\geq 0}{-r\choose u}z^u\right) \cdot \left(\sum_{v\geq 0}{r\choose v}z^v\right),\]
which implies that for $w>0$ we have
\begin{equation}\label{eq:identity-bin*bin} \sum_{u+v=w}{-r\choose u}\cdot {r\choose v} = 0.
\end{equation}
It follows that
\[ f_i\cdot(y_1+ty_2)^r = \sum_{j=i}^d\sum_{k=0}^r {-r\choose j-i}\cdot {r\choose r-k} t^{r-k+j-i}\cdot\frac{1}{y_1^{1+(j-k)}y_2^{1+(d-r-(j-k))}}.\]
Noting that the terms where $j<k$ vanish, it follows that for the remaining terms we have $w=r-k+j-i\geq r-i>0$. Setting $u=j-i$, $v=r-k$, rearranging the terms and using \eqref{eq:identity-bin*bin} gives the desired vanishing.
\end{proof}

\subsection{The proof of cases (1)--(3) of Theorem~\ref{thm:recursive-Fdr}.} We can now apply our preliminary results to establish the recursive formulas of Theorem~\ref{thm:recursive-Fdr} in all but the last case. The assertions are immediate when $r=0$, so we assume $r\geq 1$.

{\bf Case (1)} follows by combining \eqref{eq:Fdr-for-d<r} with the $T$-isomorphism $D^dU \simeq H_{d,0}$.

{\bf Case (2): $d\geq q+r-1$.} We will construct a $T$-equivariant complex
\begin{equation}\label{eq:cx-Fdr-into-FqFd-qr} 
0 \lra \mc{F}^d_r \overset{\iota}{\lra} F^qU \oo \mc{F}^{d-q}_r \overset{\pi}{\lra} \mc{F}^{d-q}_r(q) \lra 0
\end{equation}
where $\iota$ is injective, and $\pi$ is induced by the natural surjection $F^q U \oo \mc{O}_{\PP^1}\onto\mc{O}_{\PP^1}(q)$ and hence it is surjective. Since
\[ \op{rank}\left(F^qU \oo \mc{F}^{d-q}_r\right) = 2r = \op{rank}\left(\mc{F}^d_r\right) + \op{rank}\left(\mc{F}^{d-q}_r(q)\right),\]
the middle homology of \eqref{eq:cx-Fdr-into-FqFd-qr} must be $0$-dimensional. It follows from \eqref{eq:rank-det-Fdr} that the determinant of the complex is trivial, hence the complex is in fact exact. Conclusion (2) is then a consequence of the identification
\[ \ker\left(F^q U \oo \mc{O}_{\PP^1}\onto\mc{O}_{\PP^1}(q)\right) = \mc{R}^{\oo q} = L_{q,q}(-q).\]

To construct \eqref{eq:cx-Fdr-into-FqFd-qr} we use the basis $\{x_1^q,x_2^q\}$ of $F^qU$, as well as the $q$-iterates of the maps in \eqref{eq:mult-yi-Fdr}, in order to define $\iota$ via
\[ \iota(f) = x_1^q \oo y_1^q(f) + x_2^q \oo y_2^q(f)\quad\text{ for }f\text{ a local section of }\mc{F}^d_r.\]
Using the identity
\[ x_1^qy_1^q+x_2^qy_2^q = \omega^q,\]
together with the hypothesis $r\leq q$, it follows that
\[ \pi(\iota(f)) = \omega^q\cdot f = \omega^{q-r}\cdot \omega^r \cdot f = 0,\]
hence $\iota$ and $\pi$ indeed define a complex. To prove that $\iota$ is injective, we restrict to the affine chart $D_+(x_1)$ and use Lemma~\ref{lem:local-basis-Fdr}: it is enough to check the $A$-linear independence of the sections
\[ y_2^q(f_i) = \sum_{j=i}^{d-q} {-r \choose j-i}\cdot t^{j-i} \cdot\frac{1}{y_1^{j+1}y_2^{d+1-j-q}} = \frac{1}{y_1^{i+1}y_2^{d+1-i-q}} + \text{ lower order terms}.\]
The linear independence then follows from the fact that the displayed leading terms are distinct and non-zero, which uses $d-q\geq r-1\geq i$ for $i=0,\cdots,r-1$.

{\bf Case (3): $q-1\leq d\leq q+r-1$.} If $d\leq r-1$, then since $q-1\leq d$ and $r\leq q$, we must have $r=q$ and $d=q-1$, which by our convention makes $H_{d,q}=0$. The conclusion of (3) then becomes $\mc{F}^d_r = H_{q-1,0}\oo\mc{O}_{\PP^1}$ which is a special case of (1). We may therefore assume that $d\geq r$.

By an easy adaptation of the argument in \cite{rai-vdb}*{Section~4.4}, the tautological sequence \eqref{eq:ses-on-P1} gives rise after taking truncated powers to an exact complex
\[\mc{G}^{\bullet}:\quad 0 \lra T_qD^{d+q-r}U \oo \mc{O}_{\PP^1} \overset{\cdot\omega^{q-r}}{\lra} T_qD^dU(q-r) \overset{\cdot\omega^r}{\lra} T_qD^{d-r}U(q) \overset{\cdot\omega^{q-r}}{\lra} T_q D^{d-q}U(2q-r) \lra 0,\]
where $\mc{G}^0=T_qD^{d+q-r}U\oo \mc{O}_{\PP^1}$. Using \eqref{eq:def-seq-Fdr}, we get that
\[ \mc{F}^{d-r}_{q-r}(q) = \ker\left(\mc{G}^2\lra\mc{G}^3\right) = \coker\left(\mc{G}^0\lra\mc{G}^1\right).\]
Twisting by $\mc{O}_{\PP^1}(r-q)$, we can then form the following diagram with exact rows
\[
\xymatrix{
  & 0 \ar[d] & 0 \ar[d] & 0 \ar[d] & \\
0 \ar[r] & T_qD^{d+q-r}U(r-q) \ar[r] \ar[d] & T_qD^{d}U \oo \mc{O}_{\PP^1} \ar[r] \ar[d] & \mc{F}^{d-r}_{q-r}(r) \ar[r] \ar[d] & 0\\
0 \ar[r] & \mc{F}^d_r \ar[r] \ar[d] & D^{d}U \oo \mc{O}_{\PP^1} \ar[r] \ar[d] & D^{d-r}U(r) \ar[r] \ar[d] & 0\\
0 \ar[r] & D^{d-q}U \oo \mc{R}^{\oo q} \ar[r] \ar[d] & D^{d-q}U \oo F^qU\oo  \mc{O}_{\PP^1} \ar[r] \ar[d] & D^{d-q}U(q) \ar[r] \ar[d] & 0\\
  & 0  & 0  & 0  & \\
}
\]
The middle column is exact since $d<2q$ (see \cite{gao-raicu}*{(2.7)}), while the third column is exact by \eqref{eq:def-seq-Fdr}. It follows that the first column is exact as well, and moreover it is a split exact sequence, because
\[ \Ext^1\left(\mc{O}_{\PP^1}(-q),\mc{O}_{\PP^1}(r-q)\right) = H^1(\PP^1,\mc{O}_{\PP^1}(r))=0.\]
Since $\mc{R}^{\oo q}=L_{q,q}(-q)$, conclusion (3) of our theorem is then a consequence of the $T$-isomorphisms
\[ D^{d-q}U\oo L_{q,q} \simeq H_{d,q}\quad\text{and}\quad T_qD^{d+q-r}U = \left(\bw^2 U\right)^{\oo(d-r+1)}\oo \Sym^{q-2+r-d}U\simeq H_{q-1,d-r+1},\]
where the second identification follows from \cite{gao-raicu}*{(2.6)} or \cite{sun}*{Proposition~3.5}. \qed 

Case (3) of Theorem~\ref{thm:recursive-Fdr} implies the following vanishing that will be used in the proof of case (4) later on.

\begin{corollary}\label{cor:vanish-H1Fdr}
 If $r\leq d<q$ then $H^1\left(\PP^1,\mc{F}^d_r(m)\right)=0$ for $m\geq q-r-1$.
\end{corollary} 

\begin{proof} Since $\mc{F}^d_r$ decomposes as a direct sum of line bundles, and the vanishing $H^1\left(\PP^1,\mc{O}_{\PP^1}(j)\right)=0$ implies $H^1\left(\PP^1,\mc{O}_{\PP^1}(j')\right)=0$ for all $j'\geq j$, it suffices to consider the case $m=q-r-1$. Using Lemma~\ref{lem:H01-Fdr} we get
\[H^1\left(\PP^1,\mc{F}^d_r(m)\right)^{\vee} \simeq H^0\left(\PP^1,\mc{F}^{q-1}_{r}(d-r)\right)\overset{\text{Theorem}~\ref{thm:recursive-Fdr}(3)}{\simeq} H^0\left(\PP^1,\mc{O}_{\PP^1}(d-q)\right)^{\oplus r}=0,\]
where the last vanishing follows from the hypothesis $d<q$.
\end{proof}

\subsection{Action of differential operators and consequences.} In this section we explore some consequences of the fact that the local cohomology module $M$ in Section~\ref{sec:prelim} has a natural structure of a $\mc{D}_R$-module, where $\mc{D}_R$ is the ring of differential operators on $R$. Specifically we are interested in the action of differential operators
\[ \pd^{(a_1,a_2)} = \pd_1^{(a_1)}\pd_2^{(a_2)},\quad\text{ where }\pd_i^{(a)} = \frac{1}{a!}\cdot\frac{\pd^a}{\pd y_i^a},\quad\left(\text{that is, }\pd_i^{(a)}(y_i^m)={m\choose a}y_i^{m-a}\quad\text{for }m\in\bb{Z}\right)\]
and note that this action is $S$-linear. We will disregard differentiation with respect to the $x$ variables, so the notation $\pd_i$ should not be ambiguous. The following identity will be used repeatedly in what follows:
\begin{equation}\label{eq:pda1a2-act-ymon}
\pd^{(a_1,a_2)}\left(\frac{1}{y_1^{1+d_1}y_2^{1+d_2}} \right) = (-1)^{a_1+a_2}\cdot{d_1+a_1 \choose a_1}\cdot{d_2+a_2 \choose a_2}\cdot\frac{1}{y_1^{1+d_1+a_1}y_2^{1+d_2+a_2}}
\end{equation}

\begin{lemma}\label{lem:del-a1-a2-ops}
    If $a=a_1+a_2$, then the differential operators $\pd^{(a_1,a_2)}$ induce $S$- (respectively $\mc{O}_{\PP^1}$-) linear maps
    \[\pd^{(a_1,a_2)}:F^d_r \lra F^{d+a}_{r+a},\qquad \pd^{(a_1,a_2)}:\mc{F}^d_r \lra \mc{F}^{d+a}_{r+a}.\]
\end{lemma}

\begin{proof}
    We prove the $S$-module statement, which implies the sheaf statement by passing to the associated sheaves. It follows from \eqref{eq:pda1a2-act-ymon} that
    \[ \pd^{(a_1,a_2)}(M_{d,\bullet}) \subseteq M_{d+a,\bullet}.\]
    Consider now an element $f\in F^d_r$, that is, $f\in M_{d,\bullet}$ with $\omega^r\cdot f=0$. Using the fact that $\pd^{(a_1,a_2)}$ is a differential operator of order $a$, we can write
    \[ \omega^{r+a} \cdot \pd^{(a_1,a_2)} = Q \cdot \omega^r,\]
    for some differential operator $Q$. It follows that
    \[ \omega^{r+a} \cdot \pd^{(a_1,a_2)}(f) = Q \cdot \omega^r(f) = 0,\]
    hence $\pd^{(a_1,a_2)}(f)\in F^{d+a}_{r+a}$, as desired.
\end{proof}

\begin{lemma}\label{lem:Lucas-conseq}
    If $q=p^e$, $q'=p^{e-1}$, and if $v,f$ are non-negative integers satisfying $0\leq v< p$, $vq'\leq f <q$, then the binomial coefficient ${f \choose vq'}$ is non-zero (a unit) in $\kk$.
\end{lemma}

\begin{proof}
 If we write $f = uq'+f'$ with $0\leq f'<q'$, then we have $v\leq u<p$ and by Lucas' theorem it follows that
 \[ {f \choose vq'} \equiv {u\choose v} \ \text{mod }p.\]
 Since $v\leq u<p$, this implies that ${f \choose vq'}$ is non-zero in $\kk$.
\end{proof}

We obtain the following consequence, which will play an important role in the proof of Theorem~\ref{thm:recursive-Fdr}(4).

\begin{corollary}\label{cor:Fq-Da-action-Fdr}
    Suppose that $q,q'$ are as in Lemma~\ref{lem:Lucas-conseq}, and that $a,r,d$ are integers satisfying $0\leq a<p$, $aq'\leq r\leq(a+1)q'$ and $r\leq d<q$. There exists a $T$-equivariant locally split inclusion
    \[\a:F^{q'}(H_{a,0}) \oo \mc{F}^{d-aq'}_{r-aq'} \lra \mc{F}^d_r,\]
    induced by Lemma~\ref{lem:del-a1-a2-ops} and the identification of $F^{q'}(H_{a,0})$ with the $\kk$-linear span of $\pd^{(a_1q',a_2q')}$, $a_1+a_2=a$.
\end{corollary}

\begin{proof}
    Since $\{y_1,y_2\}$ is the dual basis to $\{x_1,x_2\}$, and since the derivations $\{\pd_1,\pd_2\}$ are dual to $\{y_1,y_2\}$, it follows that the $T$-action on $U$ naturally induces an isomorphism $\kk\cdot\pd^{(u,v)} \simeq L_{u,v}$. This implies that the $\kk$-linear span of $\pd^{(a_1q',a_2q')}$, $a_1+a_2=a$, is $T$-isomorphic to $F^{q'}(H_{a,0})$. Via this identification, the map
    \begin{equation}\label{eq:def-Fq'Ha0-action}
    \a:F^{q'}(H_{a,0}) \oo \mc{F}^{d-aq'}_{r-aq'} \lra \mc{F}^d_r,\quad \pd^{(a_1q',a_2q')}\oo f \lra \pd^{(a_1q',a_2q')}(f)
    \end{equation}
    is $T$-equivariant, so it suffices to check that it is a locally split inclusion. We do so on the affine chart $D_+(x_1)$, while the analogous verification on $D_+(x_2)$ follows by symmetry. To that end, we consider as in Lemma~\ref{lem:local-basis-Fdr} the $A$-basis of $\mc{F}^{d-aq'}_{r-aq'}(D_+(x_1))$ given by 
    \[ f_i' = \frac{1}{y_1^{1+i}y_2^{1+d-aq'-i}}+(\text{lower terms}),\ i=0,\cdots,r-aq'-1.\]
    It follows from \eqref{eq:pda1a2-act-ymon} that for $a_1+a_2=a$ and for $i=0,\cdots,r-aq'-1$ one has
    \[ \pd^{(a_1q',a_2q')}(f_i') = (-1)^{aq'}\cdot{i+a_1q' \choose a_1q'}\cdot{d-a_1q'-i \choose a_2q'}\cdot\frac{1}{y_1^{1+i+a_1q'}y_2^{1+d-a_1q'-i}}+(\text{lower terms}).\]
    Since $d,r<q$, it follows from Lemma~\ref{lem:Lucas-conseq} that each of the leading coefficients $(-1)^{aq'}\cdot{i+a_1q' \choose a_1q'}\cdot{d-a_1q'-i \choose a_2q'}$ is a unit in $\kk$. Moreover, since $r-aq'\leq q'$, the leading monomials $\frac{1}{y_1^{1+i+a_1q'}y_2^{1+d-a_1q'-i}}$ are also distinct as $i,a_1$ vary with $0\leq i\leq r-aq'-1$ and $0\leq a_1\leq a$. It follows that \eqref{eq:def-Fq'Ha0-action} is injective, and the sections $\pd^{(a_1q',a_2q')}(f_i')$ generate a free $A$-module which is a direct summand of $D^dU\oo A$. It follows that on $D_+(x_1)$ the composition
    \[ F^{q'}(H_{a,0}) \oo \mc{F}^{d-aq'}_{r-aq'} \lra \mc{F}^d_r \lra D^dU\oo\mc{O}_{\PP^1}\]
    is a split inclusion, which forces \eqref{eq:def-Fq'Ha0-action} to also be a split inclusion on $D_+(x_1)$.
\end{proof}

\begin{corollary}\label{cor:F_q'-vs-F_aq'}
    If $q'=p^{e-1}$, $q=p^e$, $1\leq a<p$ and $aq'\leq d<q$ then there is an isomorphism
    \[ \a':F^{q'}(H_{a-1,0}) \oo \mc{F}^{d-(a-1)q'}_{q'} \lra \mc{F}^d_{aq'}\]
    induced by differentiation as in Corollary~\ref{cor:Fq-Da-action-Fdr}. 
\end{corollary}

\begin{proof}
    We take $r=aq'$, so that $(a-1)q'\leq r\leq aq'$ is satisfied. We can then apply Corollary~\ref{cor:Fq-Da-action-Fdr} and conclude that the map
    \[\a':F^{q'}(H_{a-1,0}) \oo \mc{F}^{d-(a-1)q'}_{r-(a-1)q'} \lra \mc{F}^d_r\] 
    is a locally split inclusion. Since both the source and the target are locally free of rank $aq'$, it follows that $\a'$ is an isomorphism.
\end{proof}

\subsection{The proof of part (4) of Theorem~\ref{thm:recursive-Fdr}.} It will be useful to restrict the possible $T$-equivariant summands that appear in $\mc{F}^d_r$. We show the following (see also \eqref{eq:char-0-GLU-Fdr} for the characteristic $0$ statement). 

\begin{lemma}\label{lem:weights-Fdr} Let $d_0,r_0\in\bb{Z}_{\geq 0}\cup\{\infty\}$. If the isomorphisms in Theorem~\ref{thm:recursive-Fdr} hold for $d\leq d_0$, $r\leq r_0$, then for all such $d,r$ we have that every $T$-equivariant summand $L_{u,v}(i)$ of $\mc{F}^d_r$ satisfies $d-r+1\leq u,v\leq d$. 
\end{lemma}

\begin{proof}
    The assertion is vacuous for $r=0$. For $r\geq 1$, we argue by induction on $d$, taking $q$ to be the smallest power of $p$ with $r\leq q$ and considering the four cases of Theorem~\ref{thm:recursive-Fdr}.
    
    In case (1) we have that $L_{u,v}$ is a $T$-equivariant summand of $H_{d,0}$, hence $u,v\geq 0$ and $u+v=d$. Since $d-r+1\leq 0$, the desired conclusion follows.

    In case (2), we have that $L_{u-q,v-q}(i+q)$ is a $T$-equivariant summand of $\mc{F}^{d-q}_r$, hence by induction we have 
    \[d-q-r+1\leq u-q,v-q\leq d-q,\] 
    which implies the desired inequalities for $u,v$.

    In case (3) we have two possibilities:
    \begin{itemize}
        \item If $L_{u,v}$ is a $T$-equivariant summand of $H_{d,q}$ then $u+v=d+q$, and $u,v\geq q\geq d-r+1$. In particular, 
        \[u,v \leq (u+v)-q = d,\] 
        as desired.
        \item If $L_{u,v}$ is a $T$-equivariant summand of $H_{q-1,d-r+1}$ then $u,v\geq d-r+1$ and $u+v=d+q-r$, hence
        \[u,v \leq (u+v)-(d-r+1) = q-1 \leq d,\] 
        as desired. 
    \end{itemize}
    In case (4) we have again two possibilities:
    \begin{itemize}
     \item If $L_{u,v}(i)$ is a $T$-equivariant summand of $F^{q'}(H_{a,0}) \oo \mc{F}^{d-aq'}_{r-aq'}$, then there exist $u_0,v_0\geq 0$, $u_0+v_0=a$, such that $L_{u-u_0q',v-v_0q'}(i)$ is a $T$-equivariant summand of $\mc{F}^{d-aq'}_{r-aq'}$. The lower bound $u,v\geq d-r+1=(d-aq')-(r-aq')+1$ follows by induction. For the upper bound, we have again by induction
     \[ u = (u-u_0q') + u_0q' \leq (d-aq')+u_0q' \leq d,\]
     and similarly $v\leq d$, as desired.
     \item If $L_{u,v}(i)$ is a $T$-equivariant summand of $F^{q'}(H_{a-1,0}) \oo \mc{F}^{d-r+q'}_{(a+1)q'-r}(r-aq')$, then there exist $u_0,v_0\geq 0$, $u_0+v_0=a-1$, such that $L_{u-u_0q',v-v_0q'}(i+aq'-r)$ is a $T$-equivariant summand of $\mc{F}^{d-r+q'}_{(a+1)q'-r}$. For the upper bound, we obtain using induction and the inequality $r\geq aq'$ that
    \[ u = (u-u_0q') + u_0q' \leq (d-r+q')+u_0q' \leq d-r+aq' \leq d,\]
    and similarly $v\leq d$. For the lower bound, we have
    \[ u \geq (u-u_0q') \geq (d-r+q') - ((a+1)q'-r) + 1 = d-aq'+1\geq d-r+1,\]
    and similarly $v\geq d-r+1$, concluding our proof.\qedhere
    \end{itemize}
\end{proof}

We next identify some trivial summands of the sheaf $\mc{F}^d_r$ under the hypothesis of Corollary~\ref{cor:Fq-Da-action-Fdr}. For $d_1,d_2\geq 0$ and $q'$ as above, we write
\begin{equation}\label{eq:d12-divrem-byq}
 d_1=u_1q'+v_1, \quad d_2=u_2q'+v_2,\quad\text{with }0\leq v_1,v_2<q',
\end{equation}
and note that if $d=d_1+d_2$ satisfies $d\geq aq'$ then $u_1+u_2\geq a-1$, and the equality $u_1+u_2=a-1$ may hold when $d\leq (a+1)q'-2$.

\begin{lemma}\label{lem:triv-summand-Fdr}
 Suppose that $aq'\leq d\leq (a+1)q'-2$ for some $1\leq a<p$, and define the subspace $W\subseteq D^dU$ given by
 \[ W = \op{Span}_{\kk}\left\{\left.\frac{1}{y_1^{1+d_1}y_2^{1+d_2}} \right| d_1+d_2=d,\ u_1+u_2=a-1 \right\},\]
 where $u_1,u_2$ are determined by $d_1,d_2$ via \eqref{eq:d12-divrem-byq}. We have an isomorphism $W\simeq F^{q'}(H_{a-1,0}) \oo H_{q'-1,d-aq'+1}$ which is $T$-equivariant, and moreover
 \[ \mc{W}^d_r := W\oo\mc{O}_{\PP^1}\text{ is a direct summand of }\mc{F}^d_r\quad\text{ if }r\geq aq'.\]
\end{lemma}

\begin{proof}
 The $T$-equivariant description of $W$ follows (recalling the convention \eqref{eq:forget-bw2Udual}) from the identities
 \[ \bigoplus_{u_1+u_2=a-1} L_{u_1q',u_2q'} = F^{q'}(H_{a-1,0})\quad\text{ and }\quad \bigoplus_{\substack{0\leq v_1,v_2\leq q'-1 \\ v_1+v_2 = d-(a-1)q'}} L_{v_1,v_2} = H_{q'-1,d-aq'+1}.\]

It is clear that $W\oo\mc{O}_{\PP^1}$ is a direct summand of $D^dU\oo\mc{O}_{\PP^1}$, so to show that $\mc{W}^d_r$ is a summand of~$\mc{F}^d_r$, it suffices to check that $\omega^r$ annihilates $W$. Since
\[ \omega^r = \omega^{r-aq'} \cdot \omega^{aq'} = \omega^{r-aq'} \cdot \left(x_1^{q'}y_1^{q'} + x_2^{q'}y_2^{q'} \right)^a,\]
it follows that each term in the monomial expansion of $\omega^r$ is divisible by some $y_1^{a_1q'}y_2^{a_2q'}$ with $a_1+a_2=a$. It follows that if $u_1+u_2=a-1$ then for each such $a_1,a_2$, either $a_1>u_1$ or $a_2>u_2$, which in turn implies $d_1-a_1q'<0$ or $d_2-a_2q'<0$. We conclude that $y_1^{a_1q'}y_2^{a_2q'}\cdot W=0$, hence $\omega^r\cdot W=0$, as desired.
\end{proof}

Notice that the isomorphism in Corollary~\ref{cor:F_q'-vs-F_aq'} restricts to an isomorphism
\[ \alpha' : F^{q'}(H_{a-1,0}) \oo \mc{W}^{d-(a-1)q'}_{q'} \lra \mc{W}^d_{aq'}.\]

\begin{lemma}\label{lem:twisted-Euler}
    Let $q,q'$ be as defined in Lemma~\ref{lem:Lucas-conseq}, let $d'=d-(a-1)q'$, and assume that $(a+1)q'-1\leq d<q$. We define
    \[ E^{(q')} =\pd^{(q',0)}y_1^{q'}+\pd^{(0,q')}y_2^{q'},\]
    and consider $Q=\pd^{(a_1'q',a_2'q')}$, with $a_1'+a_2'=a-1$. For every $f\in D^{d'}U$, we have
    \[ Q(E^{(q')}(f)) = E^{(q')}(Q(f)) + (a-1)Q(f).\]
\end{lemma}

\begin{proof}
 We first note that if $d_1,d_2$ are as in \eqref{eq:d12-divrem-byq} and $d_1+d_2<q$, then
 \begin{equation}\label{eq:Eqpr-action}
  E^{(q')}\left(\frac{1}{y_1^{1+d_1}y_2^{1+d_2}}\right) \overset{\eqref{eq:pda1a2-act-ymon}}{=} (-1)^{q'}\cdot\left[{d_1\choose q'} + {d_2\choose q'} \right] \frac{1}{y_1^{1+d_1}y_2^{1+d_2}} = (-u_1-u_2)\frac{1}{y_1^{1+d_1}y_2^{1+d_2}}
 \end{equation}
where the last equality uses the fact that $(-1)^{q'}=-1$ in $\kk$, and Lucas' theorem. To prove the lemma, it suffices to consider 
\[f=\frac{1}{y_1^{1+d'_1}y_2^{1+d'_2}},\text{ where }d'_1+d'_2=d'.\] 
If we write $d'_i = u'_iq'+v'_i$, $d_i = d'_i+a'_iq'$ and use \eqref{eq:d12-divrem-byq}, we get $u_i=u'_i+a'_i$, and the hypothesis $d<q$ implies that $u_i<p$. We can then use \eqref{eq:pda1a2-act-ymon} and Lucas' theorem to obtain
\[ Q(f)= (-1)^{(a_1'+a_2')q'}{d'_1+a'_1q'\choose a'_1q'}\cdot{d'_2+a'_2q'\choose a'_2q'}\cdot \frac{1}{y_1^{1+d'_1+a'_1q'}y_2^{1+d'_2+a'_2q'}} = (-1)^{a-1}{u_1\choose a'_1}\cdot{u_2\choose a'_2}\cdot \frac{1}{y_1^{1+d_1}y_2^{1+d_2}}.\]
Combining the above with \eqref{eq:Eqpr-action} it follows that
\[ Q(E^{(q')}(f)) = (-u'_1-u'_2)\cdot Q(f) =(-1)^a(u'_1+u'_2)\cdot{u_1\choose a'_1}\cdot{u_2\choose a'_2}\cdot \frac{1}{y_1^{1+d_1}y_2^{1+d_2}}.\]
We next compute, using the same notation and reasoning as above,
\[ E^{(q')}(Q(f)) \overset{\eqref{eq:pda1a2-act-ymon}}{=} (-1)^{a-1} {u_1\choose a'_1}\cdot{u_2\choose a'_2}\cdot E^{(q')}\left[\frac{1}{y_1^{1+d_1}y_2^{1+d_2}}\right] \overset{\eqref{eq:Eqpr-action}}{=} (-1)^a(u_1+u_2)\cdot{u_1\choose a'_1}\cdot{u_2\choose a'_2}\cdot \frac{1}{y_1^{1+d_1}y_2^{1+d_2}}.\]
The desired conclusion follows from the identity 
\[(-1)^a(u'_1+u'_2) = (-1)^a(u_1+u_2) + (-1)^{a-1}(a-1),\] 
which in turn holds since $(u_1-u'_1)+(u_2-u_2') = a'_1+a'_2=a-1$.
\end{proof}

We are now ready to verify the last case of Theorem~\ref{thm:recursive-Fdr}.

\begin{proof}[Proof of Theorem~\ref{thm:recursive-Fdr}(4)] For $r=q'$ the asserted formula is tautological. We may therefore assume $r>q'$ and argue by induction on $d$. Consider first the case when $d\leq (a+1)q'-2$. The filtration \cite{gao-raicu}*{(2.7)} reduces to a short exact sequence
\[ 0 \lra F^{q'}(D^{a-1}U) \oo T_{q'}D^{d-(a-1)q'}U \lra D^d U \lra F^{q'}(D^a U) \oo D^{d-aq'}U  \lra 0.\]
Using the $T$-equivariant identifications $D^{a-1}U = H_{a-1,0}$ and $T_{q'}D^{d-(a-1)q'}U=H_{q'-1,d-aq'+1}$, it follows that $F^{q'}(D^{a-1}U) \oo T_{q'}D^{d-(a-1)q'}U=W$, where $W$ is as defined in Lemma~\ref{lem:triv-summand-Fdr}. It follows that the map $\Delta$ in \eqref{eq:def-seq-Fdr} factors through $F^{q'}(D^a U) \oo D^{d-aq'}U$, and we can write
\[ \mc{F}^d_r = \mc{W}^d_r \oplus \ol{\mc{F}}^d_r,\quad\text{where}\quad \ol{\mc{F}}^d_r = \ker\left(F^{q'}(D^a U) \oo D^{d-aq'}U\oo\mc{O}_{\PP^1} \overset{\ol{\Delta}}{\lra} D^{d-r}U(r) \right)\]
We can further factor $\ol{\Delta}$ as
\[ F^{q'}(D^a U) \oo D^{d-aq'}U\oo\mc{O}_{\PP^1} \overset{\ol{\Delta}'}{\lra} F^{q'}(D^a U) \oo D^{d-r}U(r-aq') \overset{\ol{\Delta}''}{\lra} D^{d-r}U(r), \]
where $\ol{\Delta}''$ is the surjection induced by
\[ D^a U \oo \mc{O}_{\PP^1} = \Sym^a U \oo \mc{O}_{\PP^1} \onto \mc{O}_{\PP^1}(a)\]
whose kernel is $\bw^2 U \oo \Sym^{a-1}U \oo \mc{O}_{\PP^1}(-1) = H_{a,1}(-1)$, and $\ol{\Delta}'$ is obtained by tensoring with $F^{q'}(D^a U)=F^{q'}(H_{a,0})$ from the surjection defining the short exact sequence (see \eqref{eq:def-seq-Fdr})
\[ 0 \lra \mc{F}^{d-aq'}_{r-aq'} \lra D^{d-aq'}U\oo\mc{O}_{\PP^1} \lra D^{d-r}U(r-aq') \lra 0.\]
Since $\ker(\ol{\Delta})$ is an extension of $\ker(\ol{\Delta}'')$ by $\ker(\ol{\Delta}')$, it follows that we have a short exact sequence
\begin{equation}\label{eq:ses-Fbardr}
 0 \lra F^{q'}(H_{a,0})\oo \mc{F}^{d-aq'}_{r-aq'} \lra \ol{\mc{F}}^d_r  \lra F^{q'}(H_{a,1}) \oo H_{d-r,0}(r-(a+1)q') \lra 0
\end{equation}
Using the reinterpretation of Theorem~\ref{thm:recursive-Fdr}(4) as in Corollary~\ref{cor:rec-Fdr-sometimes}(i), along with the identification $F^{q'}(H_{a,1}) \oo H_{d-r,0}=F^{q'}(H_{a-1,0}) \oo H_{d-r+q',q'}$, it suffices to verify that \eqref{eq:ses-Fbardr} splits. It is enough to show
\[ \Ext^1\left(\mc{O}_{\PP^1}(r-(a+1)q'),\mc{F}^{d-aq'}_{r-aq'}\right) = H^1\left(\PP^1,\mc{F}^{d-aq'}_{r-aq'}((a+1)q'-r)\right) = 0,\]
which follows from Corollary~\ref{cor:vanish-H1Fdr} (with $q=q'$) using the inequalities $r-aq'\leq d-aq' < q'$.

It remains to consider the case when $(a+1)q'-1\leq d\leq q-1$. We will construct below a $T$-equivariant exact sequence
\begin{equation}\label{eq:case4-Frd-seq} 
0 \lra F^{q'}(H_{a,0}) \oo \mc{F}^{d-aq'}_{r-aq'} \overset{\a}{\lra} \mc{F}^d_r \lra F^{q'}(H_{a-1,0}) \oo \mc{F}^{d-r+q'}_{(a+1)q'-r}(r-aq') \lra 0.
\end{equation}
To prove that \eqref{eq:case4-Frd-seq} is split, it will then suffice to show that there are no $T$-invariants in
\[ \Ext^1\left(F^{q'}(H_{a-1,0}) \oo \mc{F}^{d-r+q'}_{(a+1)q'-r}(r-aq'),F^{q'}(H_{a,0}) \oo \mc{F}^{d-aq'}_{r-aq'} \right).\]
Suppose otherwise, and consider $T$-equivariant summands $L_{u_1,v_1}(i_1)$ of $\mc{F}^{d-r+q'}_{(a+1)q'-r}$, and $L_{u_2,v_2}(i_2)$ of $\mc{F}^{d-aq'}_{r-aq'}$, such that there are non-zero $T$-invariants in
\[ \Ext^1\left(F^{q'}(H_{a-1,0}) \oo L_{u_1,v_1}(i_1+r-aq'),F^{q'}(H_{a,0}) \oo L_{u_2,v_2}(i_2) \right),\]
or equivalently, there are non-zero $T$-invariants in
\[ F^{q'}\left(H_{a,0}\oo H_{a-1,0}^{\vee}\right) \oo L_{u_2-u_1,v_2-v_1} \oo H^1\left(\PP^1,\mc{O}_{\PP^1}(i_2-i_1-r+aq')\right).\]
By Lemma~\ref{lem:weights-Fdr}, we have that 
\[u_2-u_1,v_2-v_1 \leq (d-aq') - \left((d-r+q')-((a+1)q'-r)+1 \right) = -1.\]
Since every $T$-equivariant summand $L_{u_3,v_3}$ of $H^1\left(\PP^1,\mc{O}_{\PP^1}(i_2-i_1-r+aq')\right)$ satisfies $u_3,v_3\leq -1$, it follows that for some $(u_3,v_3)$ there exists a summand of $F^{q'}\left(H_{a,0}\oo H_{a-1,0}^{\vee}\right)$ of the form $L_{u,v}$, where $u=u_1-u_2-u_3>0$ and $v=v_1-v_2-v_3>0$. We must then have that $q'|u,v$, and 
\[ 0 < u+v = aq'-(a-1)q' = q',\]
which is impossible.

To construct \eqref{eq:case4-Frd-seq} we let $d'=d-(a-1)q'$, $r'=r-aq'$, and build a commutative diagram
\begin{equation}\label{eq:comm-square-part4}
\begin{gathered}
\xymatrix{
 F^{q'}(H_{a-1,0}) \oo \mc{F}^{d'}_{r'} \ar[r]^(.7){\a'} \ar[d]_{\b'} & \mc{F}^d_{aq'} \ar[d]^{\b} \\
 F^{q'}(H_{a,0}) \oo \mc{F}^{d'-q'}_{r'} \ar[r]_(.7){\a} & \mc{F}^d_{r} \\
}
\end{gathered}
\end{equation}
where we think of $F^{q'}(H_{a-1,0})$ and $F^{q'}(H_{a,0})$ as spaces of differential operators as in Corollary~\ref{cor:Fq-Da-action-Fdr}. More precisely, we define $\a$ as in Corollary~\ref{cor:Fq-Da-action-Fdr} (note that $d'-q'=d-aq'$), which is a locally split inclusion because $d<q$. We define $\a'$ as the composition of the inclusion $\mc{F}^{d'}_{r'}\subseteq\mc{F}^{d'}_{q'}$ with the isomorphism from Corollary~\ref{cor:F_q'-vs-F_aq'}, and note that using \eqref{eq:ses-Fdrs} we have $\coker(\a') \simeq F^{q'}(H_{a-1,0})\oo\mc{F}^{d'-r'}_{q'-r'}(r')$. Using \eqref{eq:mult-yi-Fdr} and the inclusion $\mc{F}^d_{aq'}\subseteq\mc{F}^d_r$, we define $\beta'$, $\beta$ via
\[ \beta'(Q \oo f) = Q\cdot\pd^{(q',0)} \oo y_1^{q'}\cdot f + Q\cdot\pd^{(0,q')} \oo y_2^{q'}\cdot f,\qquad\text{and}\qquad \beta(g) = E^{(q')}(g) + (a-1)g.\]
The map $\beta$ is well-defined, since on $M_{d,\bullet}$ we have $\omega^{aq'}E^{(q')}=(E^{(q')}-a)\omega^{aq'}$. The commutativity of \eqref{eq:comm-square-part4} follows from Lemma~\ref{lem:twisted-Euler}.

The existence of the short exact sequence \eqref{eq:case4-Frd-seq} now reduces to proving that $\coker(\a)\simeq\coker(\a')$. Notice that $\coker(\a),\coker(\a')$ are locally free sheaves of the same rank, namely
\[ a(q'-r') = r - (a+1)r',\]
so it suffices to check that the natural map $\coker(\a')\lra\coker(\a)$ is surjective. This in turn reduces to proving that $\op{Im}(\a) + \op{Im}(\beta) = \mc{F}^d_r$. We check this locally on the chart $D_+(x_1)$, since the analogous statement on $D_+(x_2)$ follows by symmetry. To that end, it suffices to show that for each $k=0,\cdots,r-1$, $\op{Im}(\a) + \op{Im}(\beta)$ contains a section of $\mc{F}^d_r(D_+(x_1))$ of the form
\[ g_k = \text{unit}\cdot\frac{1}{y_1^{1+k}y_2^{1+d-k}}+(\text{lower terms}).\]
For $k = 0,\cdots,aq'-1$, we can take $g_k = \beta(f_k)$, where $\{f_0,\cdots,f_{aq'-1}\}$ is the $A$-basis of $\mc{F}^d_{aq'}(D_+(x_1))$ constructed in Lemma~\ref{lem:local-basis-Fdr} (in the special case $r=aq'$). To see that $g_k$ has the desired form, we set $d_1=k$, $d_2=d-k$ in \eqref{eq:Eqpr-action} to conclude that
\[ g_k = (a-1-u_1-u_2)\cdot\frac{1}{y_1^{1+k}y_2^{1+d-k}}+(\text{lower terms})\]
The assumption $d\geq(a+1)q'-1$ forces $u_1+u_2\geq a$, while $d<q$ implies $u_1+u_2<p$. We conclude that $-p < a-1-u_1-u_2 \leq -1$, that is, $(a-1-u_1-u_2)$ is a unit in $\kk$.

For $k=aq'+i$ with $0\leq i\leq r'-1$, we can take 
\[ g_k = \a\left(\pd^{(aq',0)}\oo f_i'\right) = \pd^{(aq',0)}(f_i') = (-1)^{aq'}{k\choose aq'}\cdot\frac{1}{y_1^{1+k}y_2^{1+d-k}}+(\text{lower terms}),\]
where $f_i'$ are as defined in the proof of Corollary~\ref{cor:Fq-Da-action-Fdr}, and ${k\choose aq'}$ is a unit by Lemma~\ref{lem:Lucas-conseq}.
\end{proof}

\section{Splitting for jet bundles}
\label{sec:split-jet}

We can now apply Theorem~\ref{thm:recursive-Fdr} and the relationship with jet bundles described in Section~\ref{subsec:pparts} in order to obtain the splitting type of $\mc{P}^k(\mc{O}_{\PP^1}(d))$ in characteristic $p>0$. We prove the following $T$-equivariant strengthening of Theorems~\ref{thm:sigma} and~\ref{thm:sigma-from-xi}.

\begin{theorem}\label{thm:Tsplitting-Pkm}
    Let $p=\op{char}(\kk)>0$, $k\geq0$, $d\in\bb{Z}$. If we define $\sigma$ as in \eqref{eq:def-sigma-from-xi} then $\sigma$ is an involution and we have a $T$-equivariant splitting
    \begin{equation}\label{eq:PkOd-T-splitting}
    \mc{P}^k(\mc{O}_{\PP^1}(d)) = \bigoplus_{i=0}^k L_{i,\sigma(i)}\oo\mc{O}_{\PP^1}(d-i-\sigma(i)).
    \end{equation}
\end{theorem}

\begin{proof}
We let $r=k+1$ and recall the notation in Section~\ref{subsec:pparts}. If $r\leq q=p^e$ then the two projections $\pi_1,\pi_2:\PP^{(r)}\lra \PP^1$ have the same composition with the $q$-power Frobenius morphism on $\PP^1$: indeed, if $s$ is a local section of $\mc{O}_{\PP^1}$ then
\[s^q\oo 1 - 1 \oo s^q = (s\oo 1 - 1 \oo s)^q \in\mc{I}^q \subseteq\mc{I}^r,\]
so $s^q\oo 1 = 1 \oo s^q$ on $\PP^{(r)}$. It follows that $\pi_2^*\mc{O}_{\PP^1}(q)=\pi_1^*\mc{O}_{\PP^1}(q)$, and the projection formula gives a natural ($T$-equivariant) identification
\[
\mc{P}^{r-1}(\mc{O}_{\PP^1}(d+q))
= \mc{P}^{r-1}(\mc{O}_{\PP^1}(d))\oo\mc{O}_{\PP^1}(q).
\]
Combining this with the periodicity $\xi_r^{d+1}=\xi_r^{d+1+q}$ in \eqref{rule2}, we may replace $d$ by any integer congruent to it modulo $q$.

We argue by induction on $r$, the case $r=1$ being immediate. For $r>1$, we take $q=p^e$ to be the smallest power of $p$ with $r\leq q$, set $q'=p^{e-1}$, and assume that $-1\leq d\leq q-2$. If $d\leq r-1$, the conclusion follows from Case~2 in Section~\ref{subsec:pparts} and \eqref{eq:basecase-xi}.

Suppose now that $r\leq d\leq q-2$, and write $r=aq'+b$, with $1\leq a<p$ and $0\leq b<q'$. Dualizing Theorem~\ref{thm:recursive-Fdr}(4) and using \eqref{eq:Pkm-vs-Fdr-case1}, we get
\begin{equation}\label{eq:Pr-1Od-inductive}
\begin{aligned}
\mc{P}^{r-1}(\mc{O}_{\PP^1}(d))\simeq{}&
F^{q'}(H_{a,0})\oo
\mc{P}^{b-1}(\mc{O}_{\PP^1}(d-aq'))\\
&\oplus L_{b,b}\oo F^{q'}(H_{a-1,0})\oo
\mc{P}^{q'-b-1}(\mc{O}_{\PP^1}(d-r+q'))\oo\mc{O}_{\PP^1}(-b),
\end{aligned}
\end{equation}
where we used the identity $F^{q'}(H_{c,0})^{\vee} \oo L_{cq',cq'} = F^{q'}(H_{c,0})$ for $c=a$ in the first summand, and for $c=a-1$ in the second summand. Since $b,q'-b<r$, we can apply the inductive hypothesis to the two smaller rank jet bundles in \eqref{eq:Pr-1Od-inductive}. Writing $i=i_0q'+i_1$, with $0\leq i_1<q'$, we explain next how the first summand in \eqref{eq:Pr-1Od-inductive} accounts for the indices $i$ in \eqref{eq:PkOd-T-splitting} with $i_1<b$, and the second for those with $i_1\geq b$. 

The $T$-equivariant line bundle summands in $F^{q'}(H_{a,0})\oo
\mc{P}^{b-1}(\mc{O}_{\PP^1}(d-aq'))$ are parametrized by pairs $(i_0,i_1)$ where $0\leq i_0\leq a$, $0\leq i_1<b$, and are given by
\begin{equation}\label{eq:summand-type1}
F^{q'}(L_{i_0,a-i_0})\oo L_{i_1,j_1}(d'-i_1-j_1),\quad\text{where }d'=d-aq',\ j_1=\xi^{d'+1}_b(i_1)+b-1-i_1.
\end{equation}
Writing $i=i_0q'+i_1$, we get using rules \eqref{rule2}, \eqref{rule3} that $\xi^{d'+1}_b(i_1)=\xi^{d+1}_r(i)$ and therefore
\[(a-i_0)q'+j_1=aq'-i_0q'+\xi^{d+1}_r(i)+b-1-i_1=\xi^{d+1}_r(i)+r-1-i=\sigma(i),\text{ and}\]
\[d'-i_1-j_1=d-(aq'+b-1+\xi^{d'+1}_b(i_1))=d-(r-1+\xi^{d+1}_r(i)) = d-i-\sigma(i).\]
It follows that the line bundle in \eqref{eq:summand-type1} is $L_{i,\sigma(i)}(d-i-\sigma(i))$.

The $T$-equivariant line bundle summands in $L_{b,b}\oo F^{q'}(H_{a-1,0})\oo
\mc{P}^{q'-b-1}(\mc{O}_{\PP^1}(d-r+q'))\oo\mc{O}_{\PP^1}(-b)$ are parametrized by pairs $(i_0,i'_1)$ where $0\leq i_0\leq a-1$, $0\leq i'_1<q'-b$, and are given by
\begin{equation}\label{eq:summand-type2}
L_{b,b}\oo F^{q'}(L_{i_0,a-1-i_0})\oo L_{i'_1,j'_1}(d-r+q'-i'_1-j'_1-b),\quad\text{where }j'_1=\xi^{d-r+q'+1}_{q'-b}(i'_1)+(q'-b-1)-i'_1.
\end{equation}
Writing $i_1=b+i'_1$ and $i=i_0q'+i_1$, we get $b\leq i_1<q'$ and $b+q'i_0+i_1'=i$. Moreover, rules \eqref{rule2}, \eqref{rule3} imply that $\xi^{d-r+q'+1}_{q'-b}(i'_1)=\xi^{d+1}_r(i)$ and therefore
\[b+q'(a-1-i_0)+j'_1 = r - 1 - i + \xi^{d+1}_r(i)=\s(i),\text{ and}\]
\[d-r+q'-i'_1-j'_1-b = d - r - (\xi^{d+1}_r(i)-1) = d-i-\sigma(i).\]
It follows that the line bundle in \eqref{eq:summand-type2} is again $L_{i,\sigma(i)}(d-i-\sigma(i))$. This implies the desired equation \eqref{eq:PkOd-T-splitting}.

Since $\mc{P}^k(\mc{O}_{\PP^1}(d))$ is $\GL(U)$-equivariant, it follows that its $T$-equivariant splitting type is invariant under the Weyl group $\mf{S}_2$, which acts by coordinate permutations (exchanging $L_{u,v}$ with $L_{v,u}$). Formula \eqref{eq:PkOd-T-splitting} then forces $\sigma$ to be an involution, concluding our proof.
\end{proof}

Theorem~\ref{thm:Tsplitting-Pkm} proves Theorems~\ref{thm:sigma} and~\ref{thm:sigma-from-xi} in positive characteristic. The characteristic-zero case of Theorem~\ref{thm:sigma} follows from Section~\ref{subsec:pparts}.

\begin{example}\label{ex:equivariant-principal-parts}
We revisit Example~\ref{ex:higher-principal-parts}. Theorem~\ref{thm:recursive-Fdr}(4), with $q=9$, $q'=3$, and $a=1$, gives
\[
\mc{F}^6_5\simeq F^3U\oo\mc{F}^3_2\oplus\mc{F}^4_1(2).
\]
Theorem~\ref{thm:recursive-Fdr}(3) applied with $q=3$ gives $\mc{F}^3_2\simeq L_{3,3}(-3)\oplus L_{2,2}(-1)$, while $\mc{F}^4_1=\mc{R}^{\oo4}=L_{4,4}(-4)$. Since $F^3U=L_{3,0}\oplus L_{0,3}$, dualizing and tensoring by $L_{6,6}$ as in \eqref{eq:Pkm-vs-Fdr-case1} yields
\[
\mc{P}^4(\mc{O}_{\PP^1}(6))\simeq
(L_{0,3}\oplus L_{3,0})(3)\oplus L_{2,2}(2)\oplus(L_{1,4}\oplus L_{4,1})(1).
\]
Forgetting the torus action, we recover the nonequivariant splitting in Example~\ref{ex:higher-principal-parts}.
\end{example}

We conclude by discussing the characterization of semistability for jet bundles.

\begin{proof}[Proof of Theorem~\ref{thm:semi-stability}]
The characteristic $p=0$ assertion follows from \eqref{eq:d-1-to-k} and Remark~\ref{rem:principal-parts}. We will therefore assume that $p>0$, and set $r=k+1$. By Theorems~\ref{thm:sigma} and~\ref{thm:sigma-from-xi}, the semistability of $\mc{P}^k(\mc{O}_{\PP^1}(d))$ is equivalent to $\xi^{d+1}_r\equiv 0$, that is $\xi_r^{d+1}(i)=0$ for all $i$. Note that $q'<r\leq q$, so by rule \eqref{rule2} we may assume $0\leq d<q$. It follows that $d=d_0q'+d_1$, and we set $s=d+1$.

If $0<s<r$, formula~\eqref{eq:basecase-xi} gives $\xi_r^s\not\equiv0$, whereas $\xi_r^r\equiv0$.

Suppose now that $r<s\leq q$, and write $r=r_0q'+r_1$, $s=s_0q'+s_1$, with $0\leq r_1,s_1<q'$. If $0<s_1<r_1$ then
\[\xi^s_r(0) \overset{\eqref{rule3}}{=} \xi^s_{r_1}(0) \overset{\eqref{rule2}}{=} \xi^{s_1}_{r_1}(0)\overset{\eqref{eq:basecase-xi}}{=} s_1-r_1\neq 0. \]
Similarly, if $r_1<s_1<q'$ then
\[\xi^s_r(r_1) \overset{\eqref{rule3}}{=} \xi^{s-r}_{q'-r_1}(0) \overset{\eqref{rule2}}{=} \xi^{s_1-r_1}_{q'-r_1}(0)\overset{\eqref{eq:basecase-xi}}{=} s_1-q'\neq 0. \]
If $s_1=0$ then it follows from \eqref{rule3} that $\xi^s_r\equiv 0$, since
\[\xi^s_{r_1} \overset{\eqref{rule2}}{=} \xi^{0}_{r_1} \overset{\eqref{eq:basecase-xi}}{\equiv} 0\quad\text{and}\quad \xi^{s-r}_{q'-r_1} \overset{\eqref{rule2}}{=} \xi^{q'-r_1}_{q'-r_1} \overset{\eqref{eq:basecase-xi}}{\equiv} 0.\]
Similarly, if $s_1=r_1$ then it follows from \eqref{rule3} that $\xi^s_r\equiv 0$, since
\[\xi^s_{r_1} \overset{\eqref{rule2}}{=} \xi^{r_1}_{r_1} \overset{\eqref{eq:basecase-xi}}{\equiv} 0\quad\text{and}\quad \xi^{s-r}_{q'-r_1} \overset{\eqref{rule2}}{=} \xi^{0}_{q'-r_1} \overset{\eqref{eq:basecase-xi}}{\equiv} 0.\]
Putting everything together, we obtain
\[
 \xi_r^s\equiv0
 \quad\Longleftrightarrow\quad
 s\geq r\ \text{ and }\ s\equiv0\text{ or }r\pmod{q'}.
\]
The congruence $s\equiv0\text{ or }r\pmod{q'}$ is equivalent to
$d_1\in\{k_1,q'-1\}$, under which the inequality $s\geq r$
becomes equivalent to $d_0\geq k_0$, concluding our proof.
\end{proof}

\section{Extremal twists for line bundle summands}
\label{sec:extremal-summands}

Although our description of the splitting type of jet bundles is fundamentally recursive, restricting our attention to the extremal twists of the line bundle summands yields more explicit consequences. The goal of this section is to discuss such consequences, as well as implications for the study of differential operators.

We start with the following observation, characterizing when $\xi^s_r$ attains its largest possible value.

\begin{lemma}\label{lem:max-xirs}
    We have $\xi^s_r(i)\leq r-1$ for all $i=0,\cdots,r-1$. If $q=p^e$ is minimal with the property that $r\leq q$ then the equality $\xi^s_r(i)=r-1$ holds if and only if $i=r-1$ and $s\equiv r-1\ (\op{mod }q)$.
\end{lemma}

\begin{proof}
    We argue by induction on $r$. If $r=1$ then $\xi^s_r(0)=0$ and $q=1$, so the conclusion is trivial. For the induction step, we assume that $r\geq 2$, so that $q\geq p$, and we let $q'=p^{e-1}$. Using \eqref{rule2}, we may assume further that $0\leq s\leq q-1$. 
    
    Suppose first that $0\leq s\leq r$. If $0\leq s\leq r-2$ then it follows from \eqref{eq:basecase-xi} that $\xi^s_r(i)\leq s<r-1$ for all $i$. If $s=r-1$ then \eqref{eq:basecase-xi} implies $\xi^s_r(r-1)=r-1$ and $\xi^s_r(i)=-1$ for $i\neq r-1$, while for $s=r$ we have $\xi^s_r(i)=0<r-1$ for all $i$. It follows that the equality $\xi^s_r(i)=r-1$ holds if and only if $i=s=r-1$.

    We are left considering the case $r<s<q$, for which we apply \eqref{rule3}. Since $q=p^e$ is minimal, it follows that $q'<r$, hence $r_1,q'-r_1<r$. By induction, we have $\xi^s_{r_1}(i_1)\leq r_1-1<r-1$, and $\xi^{s-r}_{q'-r_1}(i_1-r_1)\leq q'-r_1-1<r-1$. This shows that $\xi^s_r(i)<r-1$ for all $i$, concluding our proof. 
\end{proof}

Taking $r=k+1$ and $s=d+1$ we can then obtain the characterization of the splitting of \eqref{eq:can-ses-Pk-Pk-1} as follows. 

\begin{proof}[Proof of Theorem~\ref{thm:split-ses-pparts}]
    As explained in the Introduction, the splitting of \eqref{eq:can-ses-Pk-Pk-1} is equivalent to the condition $\sigma(k)=k$, which in light of \eqref{eq:def-sigma-from-xi} is further equivalent to $\xi^{d+1}_{k+1}(k)=k$. Applying Lemma~\ref{lem:max-xirs} with $r=k+1$ and $s=d+1$ yields the desired conclusion.
\end{proof}

We next explain the characterization of the extremal twists for line bundle direct summands of $\mc{P}^k(\mc{O}_{\PP^1}(d))$.

\begin{proof}[Proof of Theorem~\ref{thm:extr-summands-Pkd}] Recalling the notation \eqref{eq:def-m-pl-min} and \eqref{eq:def-D+-}, we write $D_{\pm}=D_{\pm}(k,d)$ and $m_{\pm}=m_{\pm}(k,d)$. If $d\in\mc{SS}_k$ then $d=D_-=D_+$ and $m_-=m_+=d-k$, so the desired conclusion follows. We will therefore assume that $d\not\in\mc{SS}_k$, and in particular $k\geq 1$ (since $\mc{SS}_0=\bb{Z}$). If $0\leq d< k$, then \eqref{eq:d-1-to-k} implies $m_-=-k-1$ and $m_+=0$, while Theorem~\ref{thm:semi-stability} yields $D_-=-1$ and $D_+=k$. We get $m_{\pm}=D_{\pm}-k$, as desired. This proves the theorem in characteristic zero, since in that case $\mc{SS}_k=\{d\in\bb{Z}|d<0\text{ or }d\geq k\}$.

In characteristic $p>0$, we consider $q=p^e$, where $e\geq 0$ is minimal with the property that $q>k$. It follows from Theorem~\ref{thm:semi-stability} that $\mc{SS}_k$ is invariant under translation by $q$, hence $D_{\pm}(k,d+q)=D_{\pm}(k,d)+q$. Moreover, Theorem~\ref{thm:sigma-from-xi} and \eqref{rule2} imply that $m_{\pm}(k,d+q)=m_{\pm}(k,d)+q$, so we may assume $0\leq d\leq q-1$. We have already discussed the case $0\leq d<k$, and $d=k\in\mc{SS}_k$, so it remains to consider $k<d<q$.

We let $q'=p^{e-1}$, and write $k=k_0q'+k_1$ with $0\leq k_0<p$, $-1\leq k_1<q'-1$. Since $d\not\in\SS_k$ we have that $D_-<d<D_+$. It follows from Theorem~\ref{thm:semi-stability} that either
\begin{itemize}
    \item[(i)] $D_-=aq'+k_1$ and $D_+=aq'+(q'-1)$, where $k_0\leq a\leq p-1$, or
    \item[(ii)] $D_-=aq'+(q'-1)$ and $D_+=(a+1)q'+k_1$, where $k_0\leq a\leq p-2$.
\end{itemize}

We first verify that $\mc{O}_{\PP^1}(D_--k)$ and $\mc{O}_{\PP^1}(D_+-k)$ are summands of $\mc{P}^k(\mc{O}_{\PP^1}(d))$. To do so, by Theorem~\ref{thm:sigma-from-xi} it suffices to find indices $i_{\pm}$ with $\xi^{d+1}_{k+1}(i_{\pm})=d-D_{\pm}$. In case (i), we can write $d=aq'+d_1$, where $0\leq d_1<q'$. If we take $i_-=d_1+1$ then
\[\xi^{d+1}_{k+1}(i_-)\overset{\eqref{rule3}}{=}\xi^{d-k}_{q'-1-k_1}(d_1-k_1)\overset{\eqref{rule2}}{=}\xi^{d_1-k_1}_{q'-1-k_1}(d_1-k_1)\overset{\eqref{rule1}}{=}d_1-k_1=d-D_-.\]
If we take $i_+=k_1+1$ then
\[\xi^{d+1}_{k+1}(i_+)\overset{\eqref{rule3}}{=}\xi^{d-k}_{q'-1-k_1}(0)\overset{\eqref{rule2}}{=}\xi^{d_1-k_1}_{q'-1-k_1}(0)\overset{\eqref{rule1}}{=}d_1+1-q'=d-D_+.\]
In case (ii) we have $0\leq k_1\leq q'-2$, and we can write $d=(a+1)q'+d_1$, where $0\leq d_1<k_1$. If we take $i_-=d_1+1$ then
\[\xi^{d+1}_{k+1}(i_-)\overset{\eqref{rule3}}{=}\xi^{d+1}_{k_1+1}(d_1+1)\overset{\eqref{rule2}}{=}\xi^{d_1+1}_{k_1+1}(d_1+1)\overset{\eqref{rule1}}{=}d_1+1=d-D_-.\]
If we take $i_+=0$ then
\[\xi^{d+1}_{k+1}(i_+)\overset{\eqref{rule3}}{=}\xi^{d+1}_{k_1+1}(0)\overset{\eqref{rule2}}{=}\xi^{d_1+1}_{k_1+1}(0)\overset{\eqref{rule1}}{=}d_1-k_1=d-D_+.\]

To conclude, it remains to show that every summand $\mc{O}_{\PP^1}(m)$ satisfies $D_--k\leq m\leq D_+-k$. Since $d\geq k+1$, this is equivalent by \eqref{eq:Pkm-vs-Fdr-case1} to bounding the summands $\mc{O}_{\PP^1}(-m)$ that occur in $\mc{F}^d_{k+1}$.

To prove that $m\geq D_--k$, we have to verify $H^0\left(\PP^1,\mc{F}^{d}_{k+1}(D_--k-1)\right)=0$, which by Lemma~\ref{lem:H01-Fdr} is equivalent to $H^1\left(\PP^1,\mc{F}^{D_-}_{k+1}(d-k-1)\right)=0$. We have
\[\mc{F}^{D_-}_{k+1}(d-k-1)\overset{\eqref{eq:Pkm-vs-Fdr-case1}}{\simeq}\left(\mc{P}^k(\mc{O}_{\PP^1}(D_-)) \right)^{\vee} \oo \mc{O}_{\PP^1}(d-k-1) \simeq \mc{O}_{\PP^1}(d-D_--1)^{\oplus (k+1)},\]
where the second isomorphism follows from the fact that $D_-\in\SS_k$. Since $d-D_--1\geq -1$, we obtain $H^1\left(\PP^1,\mc{F}^{D_-}_{k+1}(d-k-1)\right)=0$, as desired. 

To prove that $m\leq D_+-k$, we have to verify $H^1\left(\PP^1,\mc{F}^{d}_{k+1}(D_+-k-1)\right)=0$, which by Lemma~\ref{lem:H01-Fdr} is equivalent to $H^0\left(\PP^1,\mc{F}^{D_+}_{k+1}(d-k-1)\right)=0$. We have
\[\mc{F}^{D_+}_{k+1}(d-k-1)\overset{\eqref{eq:Pkm-vs-Fdr-case1}}{\simeq}\left(\mc{P}^k(\mc{O}_{\PP^1}(D_+)) \right)^{\vee} \oo \mc{O}_{\PP^1}(d-k-1) \simeq \mc{O}_{\PP^1}(d-D_+-1)^{\oplus (k+1)},\]
since $D_+\in\SS_k$. Since $d-D_+-1\leq -1$, we obtain $H^0\left(\PP^1,\mc{F}^{D_+}_{k+1}(d-k-1)\right)=0$, concluding our proof.
\end{proof}

The extremal twists translate directly into criteria for global generation and for the existence of global differential operators, as follows.

\begin{proof}[Proof of Corollary~\ref{cor:diffops}]
    If we decompose $\mc{P}^k(\mc{O}_{\PP^1}(d))$ as in \eqref{eq:def-m-pl-min}, and use \eqref{eq:diff-Dkde}, then it follows that
    \[\mc{D}^k_{d,e}\simeq\bigoplus_{i=0}^k\mc{O}_{\PP^1}(e-a_i).\]
    Therefore $\mc{D}^k_{d,e}$ is globally generated if and only if $e\geq a_i$ for all $i$, and it has non-zero global sections if and only if $e\geq a_i$ for some $i$, proving conclusions (1) and (2).

    Taking $e=d$ and applying Theorem~\ref{thm:sigma}, we get
    \[\mc{D}^k_{d,d}\simeq\bigoplus_{i=0}^k\mc{O}_{\PP^1}(i+\sigma(i)),\]
    which is globally generated since $i+\sigma(i)\geq 0$. Using the fact that $\dim_{\kk}H^0(\PP^1,\mc{O}_{\PP^1}(m))=m+1$ for $m\geq 0$, and that $\sigma$ is a permutation of $\{0,\cdots,k\}$, we obtain by taking global sections of $\mc{D}^k_{d,d}$ that
    \[\dim_{\kk}\op{Diff}^k(\mc{O}_{\PP^1}(d),\mc{O}_{\PP^1}(d))=\sum_{i=0}^k(i+\sigma(i)+1) = (k+1)^2.\]
    This verifies (3) and concludes our proof.
\end{proof}

\section{Cohomology of line bundles on infinitesimal thickenings}
\label{sec:coh-thickenings}

The goal of this section is to provide explicit formulas for sheaf cohomology groups of line bundles on infinitesimal thickenings of the diagonal in $\PP^1\times \PP^1$, or equivalently, of the incidence correspondence in $\PP^1\times(\PP^1)^{\vee}$. To reconcile the two points of view, we note that the tautological exact sequence on $(\PP^1)^{\vee}=\op{Proj}(\Sym(U^{\vee}))$ is obtained by dualizing \eqref{eq:ses-on-P1}, so via the natural identification of $(\PP^1)^{\vee}$ with $\PP^1$, the universal quotient bundle $\mc{O}_{(\PP^1)^{\vee}}(1)$ corresponds to $\mc{R}^{\vee}=\bw^2 U^{\vee} \oo \mc{O}_{\PP^1}(1)$. This is abstractly isomorphic to the universal quotient bundle on $\PP^1$, the only difference stems from the equivariant structure where they differ by a determinant factor. In particular, we get a natural isomorphism $\PP^{(r)}\simeq X^{(r)}$ under which $\mc{O}_{\PP^{(r)}}(a,b)$ corresponds to $L_{b,b}\oo\mc{O}_{X^{(r)}}(a,b)$. Therefore the formulas for cohomology characters that we describe here will be explained for $\PP^{(r)}$, and we leave the easy translation to $X^{(r)}$ to the interested reader. Dimension formulas for the cohomology of $\mc{O}_{\PP^{(r)}}(a,b)$ will be the same as those for $\mc{O}_{X^{(r)}}(a,b)$.

\begin{theorem}\label{thm:coh-characters-Pr}
    Let $r\geq1$, $a,b\in\bb{Z}$, and let $\sigma$ be the involution described in the Introduction for $k=r-1$ and $d=b$. Using notation \eqref{eq:def-Hab}, we have $T$-equivariant isomorphisms
    \[
    H^0\left(\PP^{(r)},\mc{O}_{\PP^{(r)}}(a,b)\right)=\bigoplus_{i+\sigma(i)\leq a+b} L_{i,\sigma(i)} \oo H_{a+b-i-\sigma(i),0},
    \]
    \[
    H^1\left(\PP^{(r)},\mc{O}_{\PP^{(r)}}(a,b)\right)=\bigoplus_{i+\sigma(i)\geq a+b+2} L_{i,\sigma(i)} \oo H_{-1,a+b+1-i-\sigma(i)}.
    \]
    In particular, $\mc{O}_{\PP^{(r)}}(a,b)$ has vanishing cohomology if and only if $i+\sigma(i)=a+b+1$ for all $i=0,\cdots,r-1$.
\end{theorem}

\begin{proof}
    If we write $\pi_1:\PP^{(r)}\lra\PP^1$ as in Section~\ref{subsec:pparts}, then $\mc{O}_{\PP^{(r)}}(a,0)=\pi_1^*(\mc{O}_{\PP^1}(a))$ and it follows from the projection formula that
    \begin{equation}\label{eq:pi1-OPr-ab}
    \pi_{1*}\left(\mc{O}_{\PP^{(r)}}(a,b)\right)=\mc{P}^{r-1}(\mc{O}_{\PP^1}(b)) \oo \mc{O}_{\PP^1}(a) \overset{\eqref{eq:PkOd-T-splitting}}{=}\bigoplus_{i=0}^{r-1} L_{i,\sigma(i)}\oo\mc{O}_{\PP^1}(a+b-i-\sigma(i)).
    \end{equation}
    For line bundles on $\PP^1$ the non-vanishing cohomology is described $T$-equivariantly by
    \[H^0\left(\PP^1,\mc{O}_{\PP^1}(d)\right)=H_{d,0}\text{ for }d\geq 0,\quad H^1\left(\PP^1,\mc{O}_{\PP^1}(d)\right)=H_{-1,d+1}\text{ for }d\leq -2.\]
    The formulas for the cohomology of $\mc{O}_{\PP^{(r)}}(a,b)$ follow using \eqref{eq:pi1-OPr-ab} from the fact that $\pi_1$ is finite, which implies
    \begin{equation}\label{eq:Hjpfwd=HjonXr}
    H^j\left(\PP^{(r)},\mc{O}_{\PP^{(r)}}(a,b)\right) = H^j\left(\PP^1,\pi_{1*}\left(\mc{O}_{\PP^{(r)}}(a,b)\right)\right)
    \end{equation}
    The final cohomology vanishing statement follows from our formulas, since there are no summands of either $H^0\left(\PP^{(r)},\mc{O}_{\PP^{(r)}}(a,b)\right)$ or $H^1\left(\PP^{(r)},\mc{O}_{\PP^{(r)}}(a,b)\right)$ with $i+\sigma(i)=a+b+1$.
\end{proof}

We can now specialize our result to the non-equivariant setting, and to line bundles with $a+b=r-2$ to obtain Theorem~\ref{thm:coho-consequences}.

\begin{proof}[Proof of Theorem~\ref{thm:coho-consequences}]
    Using the last part of Theorem~\ref{thm:coh-characters-Pr}, cohomology vanishing is equivalent to $i+\sigma(i)=a+b+1$ for all $i$, which in light of Theorem~\ref{thm:sigma} is further equivalent to $\mc{P}^{r-1}(\mc{O}_{\PP^1}(b))\simeq\mc{O}_{\PP^1}(-a-1)^{\oplus r}$. This occurs if and only if $\mc{P}^{r-1}(\mc{O}_{\PP^1}(b))$ is semistable and $-a-1=b-(r-1)$, proving the first part of Theorem~\ref{thm:coho-consequences}.

    For the second part, we write $\chi=h^0-h^1$ for the Euler characteristic (of a sheaf on a $1$-dimensional projective scheme). We have for all $d\in\bb{Z}$
    \begin{equation}\label{eq:hi-chi-Od}
    h^0\left(\PP^1,\mc{O}_{\PP^1}(d)\right)=\max(d+1,0),\quad h^1\left(\PP^1,\mc{O}_{\PP^1}(d)\right)=\max(-d-1,0),\quad \chi\left(\PP^1,\mc{O}_{\PP^1}(d)\right)=d+1.
    \end{equation}
    It follows from Theorem~\ref{thm:sigma}, \eqref{eq:pi1-OPr-ab} and \eqref{eq:Hjpfwd=HjonXr} that
    \[\chi\left(X^{(r)},\mc{O}_{X^{(r)}}(a,b)\right)=\sum_{i=0}^{r-1}(a+b+1-i-\sigma(i))=r(a+b-r+2),\]
    which implies the equality $h^0\left(X^{(r)},\mc{O}_{X^{(r)}}(a,b)\right)=h^1\left(X^{(r)},\mc{O}_{X^{(r)}}(a,b)\right)$ when $a+b=r-2$. It follows from \eqref{eq:hi-chi-Od} that
    \[h^0\left(\PP^1,\mc{O}_{\PP^1}(d)\right)+h^1\left(\PP^1,\mc{O}_{\PP^1}(d)\right)=|d+1|,\]
    which then implies that for $a+b=r-2$ we have
    \[
    \begin{aligned}
        h^0\left(X^{(r)},\mc{O}_{X^{(r)}}(a,b)\right)&=\frac{h^0\left(X^{(r)},\mc{O}_{X^{(r)}}(a,b)\right)+h^1\left(X^{(r)},\mc{O}_{X^{(r)}}(a,b)\right)}{2} \\
        &=\frac{1}{2}\cdot\sum_{i=0}^{r-1}|a+b+1-i-\sigma(i)| = \frac{1}{2}\sum_{i=0}^{r-1}|\xi^{b+1}_r(i)|,
    \end{aligned}\]
    concluding the proof.
\end{proof}

\section*{Acknowledgements}
Experiments with Macaulay2 \cite{GS} and discussions with ChatGPT \cite{ChatGPT} have provided many valuable insights. Marangone gratefully acknowledges that this research was supported in part by the Pacific Institute for the Mathematical Sciences. Marangone acknowledges the European Union’s Horizon Europe research and innovation program under the Marie Skłodowska-Curie Funding Program (Project PrIMes, Grant Agreement No. 101277460). Raicu and Reed acknowledge the support of the National Science Foundation Grant DMS-2302341. Reed also acknowledges support by the National Natural Science Foundation of China (No. 12288201) and by the National Key R$\&$D Program of China (2023YFA1009402). Part of the material in this paper is based upon work supported by the National Science Foundation under Grant No. DMS-1928930 and by the Alfred P. Sloan Foundation under grant G-2021-16778, while Raicu and Reed were in residence at the Simons Laufer Mathematical Sciences Institute (formerly MSRI) in Berkeley, California, during the Spring 2024 semester.

\end{document}